\documentclass[a4paper,11pt]{article}
\usepackage[a4paper, total={6.45in, 9in}]{geometry}
\usepackage[utf8]{inputenc} 
\usepackage[T1]{fontenc} 
\usepackage{amsmath}
\usepackage{circuitikz}
\usepackage{amsthm}
\usepackage{verbatim}
\usepackage{mathtools} 
\usepackage{extarrows} 
\usepackage{amsthm}
\usepackage{xparse}
\usepackage{hyperref}
\usepackage{cleveref}
\usepackage{setspace}
\usepackage{mathpazo} 
\usepackage{amssymb}
\usepackage{tikz-cd}
\usepackage{float}
\usepackage{graphicx}
\usepackage{tikz}
\usepackage[T1]{fontenc}
\usepackage{imakeidx}

\usepackage{xcolor}
\usepackage{color}
\usepackage{amssymb,latexsym} 
 
\numberwithin{equation}{section} 
 
\newtheorem{theorem}{Theorem}[subsection] 
\newtheorem*{teo*}{Theorem}
\newtheorem{proposition}[theorem]{Proposition} 
 
\newtheorem{corollary}[theorem]{Corollary} 
\newtheorem{lemma}[theorem]{Lemma} 
 
\theoremstyle{definition} 
 
\newtheorem{remark}[theorem]{Remark} 
\newtheorem{examples}[theorem]{Examples}
\newtheorem{example}[theorem]{Example} 
\newtheorem{definition}[theorem]{Definition} 

\usepackage[maxbibnames=99,style=alphabetic]{biblatex}
\usepackage{ytableau}
\usepackage{hyperref}
\usepackage[autostyle=true]{csquotes}

\makeindex

\begin{document}

\begin{center}
\textbf{\large A MULTIPLICATION FORMULA FOR CLUSTER CHARACTERS IN GENTLE ALGEBRAS} 
\vspace{5mm}
 \\Azzurra Ciliberti\footnote{Fakultät für Mathematik - Ruhr-Universität Bochum.\\
 \textit{Email address}: \textbf{azzurra.ciliberti@ruhr-uni-bochum.de}}

 \end{center}
 \begin{abstract}
  \noindent 
  We prove a multiplication formula for cluster characters induced by generating extensions in a gentle algebra $A$, generalizing a result of Cerulli Irelli, Esposito, Franzen, Reineke. In the case where $A$ is the gentle algebra of a triangulation $T$ of an unpunctured marked surface, this provides a representation-theoretic interpretation of the exchange relations in the cluster algebra with principal coefficients in $T$. As an application, we interpret a formula that relates cluster variables of type $B$ to cluster variables of type $A$ in the symmetric module category of the algebras arising from special triangulations of a regular polygon. 

\end{abstract}
\tableofcontents

\section*{Introduction}
\addcontentsline{toc}{section}{Introduction}

Cerulli Irelli, Esposito, Franzen, Reineke, in \cite{CEHR}, study a decomposition of the quiver Grassmannian induced by short exact sequences in the module category of the path algebra $A=kQ$, where $Q$ is an acyclic quiver and $k$ is an algebraically closed field. In particular, given two $A$-modules $X,S$ such that $\text{dim}\mathrm{Ext}^1(S,X)=1$, and a generating extension 
\begin{center}
    $\xi : 0\to X \to Y \to S \to 0$,
\end{center}
they prove that, for any $\bold{e} \in \mathbb{Z}_{\geq 0}^{Q_0}$,
\begin{equation}\label{eq_intro}
        \chi(Gr_{\bold{e}}(X \oplus S))=\chi(Gr_{\bold{e}}(Y))+ \chi(Gr_{\bold{e}-\textbf{dim}\underline{S}}(\overline{X} \oplus S/\underline{S})),
    \end{equation}
where $\chi$ is the Euler-Poincar\'e characteristic, $\overline{X}$ is the maximal submodule of $X$ such that the push-out sequence of $\xi$ does not split, and dually $\underline{S}$ is the minimal submodule of $S$ such that the pullback sequence of $\xi$ does not split. This allows them to obtain a multiplication formula between the corresponding cluster characters, and a representation-theoretic interpretation of the exchange relations in the cluster algebra associated with $Q$.

\vspace{0.3cm}

The main aim of this paper is to generalize their result to the case where $A = kQ/I$ is a gentle algebra. Gentle algebras, introduced by Assem and Skowroński in \cite{AS}, are a particularly interesting and nice class of algebras. They are string algebras, so their module categories are combinatorially described in terms of strings and bands \cite{BR,C-B}. They are tiling algebras associated with dissections of unpunctured surfaces with boundary \cite{BCS,OPS}. They are Gorenstein, tame, derived-tame, and closed under derived equivalence. Moreover, as shown in \cite{ABCP}, Jacobian algebras arising from triangulations of unpunctured surfaces with marked points are gentle.

\hspace{0.3cm}

The definitions of $\overline{X}$ and $\underline{S}$, as well as the proof of \ref{eq_intro}, heavily rely on the fact that the algebra $kQ$ is hereditary, meaning $\mathrm{Ext}_{kQ}^i(-,-) = 0$ for any $i > 1$. To generalize \ref{eq_intro} in this new setting, the idea is to define $\overline{X}$ and $\underline{S}$ differently, as follows:
\begin{center} $\overline{X} := \text{ker}(f) \subset X$; \hspace{2cm}$\underline{S} := \text{im}(g) \subseteq S$ \end{center}
where $f: X \to \tau S$ is the non-zero morphism from $X$ to $\tau S$ that does not factor through an injective $A$-module, and $g: \tau^{-1} X \to S$ is the non-zero morphism from $\tau^{-1} X$ to $S$ that does not factor through a projective $A$-module. By the Auslander-Reiten formulas (cf. \cite{libro_blu}, IV.2.13), $\overline{X}$ and $\underline{S}$ are well-defined $A$-modules. Moreover, according to \cite[Lemma 31]{CEHR}, the new definitions coincide with the original ones when $A$ is hereditary. 

\vspace{0.3cm}

With these new definitions, we obtain the first main result in this work:

\begin{teo*}\ref{dec_qg_gent_alg}
    Let $A=kQ/I$ be a gentle algebra.  Let $X$, $S$ be $A$-modules such that $\text{dim}\mathrm{Ext}^1(S,X)=1$. Let $\xi: 0 \to X \to Y \to S \to 0 \in \mathrm{Ext}^1(S,X)$ be a generating extension with a middle term $Y$. Then for any $\bold{e}\in \mathbb{Z}_{\geq0}^{Q_0}$
    \begin{equation*}
        \chi(Gr_{\bold{e}}^A(X \oplus S))=\chi(Gr_{\bold{e}}^A(Y))+ \chi(Gr_{\bold{e}-\textbf{dim}\underline{S}}^{B}(M)),
    \end{equation*} 
    where $M$ is the $\leq_{\mathrm{Ext}}$-minimum extension between $S/\underline{S}$ and $\overline{X}$ in a finite-dimensional algebra $B\supseteq A$. In particular, if $A$ is the gentle algebra of a triangulation of an unpunctured surface with marked points, then $B=A$. Moreover, $M <_{\mathrm{Ext}}\overline{X}\oplus S/\underline{S}$ if and only if $\mathrm{Ext}^1(\underline{S},X/\overline{X})=0$. 
\end{teo*}

If $A$ is hereditary, that is, $Q$ has no oriented cycles and $I=0$, by \cite[Lemma 31]{CEHR} $\mathrm{Ext}^1(\underline{S},X/\overline{X})\neq 0$ for any $X,S$ such that $\text{dim}\mathrm{Ext}^1(S,X)=1$, so we recover \ref{eq_intro}.

With the same hypothesis of Theorem \ref{dec_qg_gent_alg}, we deduce the following multiplication formula between cluster characters (up to some rescaling by $\bold{x}$-monomials, cf. Corollary \ref{cor1} for the precise statement):
\begin{equation}\label{eq_intro2}
        CC(X)CC(S)= CC(Y) + \bold{y}^{\textbf{dim} \underline{S}}CC(M).
    \end{equation} 

Furthermore, we prove the following theorem which is the second main result of this note:

\begin{teo*}[\ref{thm_ca}]
    If $A$ is the gentle algebra of a triangulation $T$ of an unpunctured surface with marked points, and in addition $X$ and $S$ are rigid and indecomposable, then \ref{eq_intro2} is an exchange relation between the cluster variables $CC(X)$ and $CC(S)$ in the cluster algebra with principal coefficients in the initial seed whose exchange matrix corresponds to $T$. 
\end{teo*}

Skein relations in cluster algebras from surfaces are investigated in other works using different techniques. For example, in \cite{CSA,MW} they come from smoothing of crossings between generalized arcs in the surface. For the punctured case, see \cite{BKK,banaian2024skeinrelationspuncturedsurfaces}. 

\vspace{0.3cm}

Finally, we apply the results of Theorem \ref{dec_qg_gent_alg} and Theorem \ref{thm_ca} to cluster algebras of type $B$, in order to    
provide a categorical interpretation of the formula that relates cluster variables of type $B_n$ to cluster variables of type $A_n$,  presented in \cite[Theorem 3.7]{ciliberti2024} (cf. Section \ref{section_B}, Theorem \ref{t1}).

Let $\mathbf{P}_{2n+2}$ be the regular polygon with $2n+2$ vertices, and let $d$ be a diameter. Let $\rho$ be the reflection of $\mathbf{P}_{2n+2}$ along $d$, and let $T'$ be a $\rho$-invariant triangulation of the polygon such that $d \in T'$. The algebra $A'=Q(T')/I(T')$ naturally associated with $T'$ is a symmetric algebra in the sense of \cite{boos2021degenerations}. On the other hand, associated with $T'$ is a cluster algebra $\mathcal{A}_\bullet^B$ of type $B_n$ with principal coefficients in the initial seed corresponding to $T'$. In \cite{ciliberti2024, 5}, the author establishes a correspondence between the cluster variables of $\mathcal{A}_\bullet^B$ and the orthogonal indecomposable $A'$-modules. 

For an orthogonal indecomposable $A'$-module $N$, $F_N$ and $\bold{g}_N$ denote the $F$-polynomial and the $\bold{g}$-vector of the cluster variable $x_N$ of $\mathcal{A}_\bullet^B$ that corresponds to $N$. On the other hand, $F_{\text{Res}(N)}$ and $\bold{g}_{\text{Res}(N)}$ are the $F$-polynomial and the $\bold{g}$-vector of the module $\text{Res}(N)$, obtained from $N$ by assigning the trivial vector space to all vertices with indices greater than $n$ (cf. Definition \ref{def_res}). The following result, obtained in \cite[Theorem 4.28]{ciliberti2024} only for the case where $A'$ is hereditary, is now proved in full generality:
\begin{teo*}\ref{cat_interpr}
Let $N$ be an orthogonal indecomposable $A'$-module. Let $D=\text{diag}(1,\dots,1,2)\in \mathbb{Z}^{n\times n}$.
\begin{itemize}
    \item [(i)] If $\text{Res}(N)=(V_i,\phi_a)$ is indecomposable as ordinary module, then 
    \begin{equation}\label{e_1_intro}
        \text{$F_N=F_{\text{Res}(N)}$,}
    \end{equation}
    and
    \begin{equation}\label{e_2_intro}
 \bold{g}_{N}=\begin{cases}
          \text{$D \bold{g}_{\text{Res}(N)}$ \hspace{1.8cm}if $\textbf{dim} V_n =0$;}\\
          \text{$D \bold{g}_{\text{Res}(N)}+\bold{e}_n$ \hspace{1cm}if $\textbf{dim} V_n \neq 0$.}
      \end{cases}
      \end{equation}
\item [(ii)]      Otherwise, $N=L\oplus \nabla L$ with $\text{dim Ext}^1(\nabla L, L)=1$, and there exists a generating extension
    \begin{center}
        $0 \to L \to G_1 \oplus G_2 \to \nabla L \to 0 $,
    \end{center}
    where $G_1$ and $G_2$ are orthogonal indecomposable $A'$-modules of type I. Then
    \begin{equation}\label{e_3_intro}
        F_N= F_{\text{Res}(N)} - \bold{y}^{\text{Res}(\textbf{dim}\underline{\nabla L})}F_{\text{Res}(M)},
    \end{equation}
    and
    \begin{equation}\label{e_4_intro}
        \bold{g}_{N}=D(\bold{g}_{\text{Res}(N)}+\bold{e}_n),
    \end{equation}
where $M$ is the $\leq_{\mathrm{Ext}}$-minimum extension in $A'$ between $\nabla L/\underline{\nabla L}$ and $\overline{L}$. Moreover, $M$ is an orthogonal $A'$-module.
 \end{itemize}  
\end{teo*}


\vspace{0.3cm}
 
The paper is structured as follows. In Section \ref{s1}, we give an overview of some basic concepts about gentle algebras and their combinatorics. In Section \ref{section_cc}, we recall the definitions of cluster character, $F$-polynomial, and $\bold{g}$-vector of a module over a finite-dimensional algebra. In Section \ref{s_acyc}, we describe the multiplication formula of \cite{CEHR} for acyclic quivers. Then, in Section \ref{s_gent_alg_cc}, we present and prove the generalization of this formula to gentle algebras. Finally, in Section \ref{section_B}, we explain how the results of Section \ref{s_gent_alg_cc} allow us to obtain an interpretation of the formula relating cluster variables of type $B_n$ and $A_n$ in the symmetric module category of the algebras associated with $\rho$-invariant triangulations of polygons.

\section{Basic notions on gentle algebras and their module categories}\label{s1}

Let $k$ be an algebrically closed field. Let $Q=(Q_0,Q_1, s,t)$ be a quiver, where $Q_0$ is the finite set of vertices, $Q_1$ is the finite set of arrows, and $s,t:Q_1 \to Q_0$ are two functions that provide the orientation $a:s(a)\to t(a)$ of the arrows from \emph{source} to \emph{target}. The \emph{path algebra} $kQ$ of $Q$ is the $k$-vector space generated by the set of all paths in $Q$. The multiplication of two paths $\alpha=a_1\cdots a_s$ and $\beta=b_1\cdots b_t$ is defined as the concatenation $\alpha\beta$ if $t(a_s)=s(b_1)$, 0 otherwise. Let $R$ be the two-sided ideal generated by the arrows of $Q$. Let $I \subseteq kQ$ be an ideal such that there exists an integer $m\geq 2$ such that $R^m \subseteq I \subseteq R^2$. Then the finite-dimensional quotient algebra $A=kQ/I$ is called \emph{quiver algebra}. We denote by $\mathrm{mod}A$ the category of finitely generated right $A$-modules. 

\begin{definition}[Gentle algebra]
    A finite dimensional algebra $A$ is \emph{gentle} if it admits a presentation $A=kQ/I$ satisfying the following conditions:
    \begin{itemize}
        \item [(G1)] Each vertex of $Q$ is the source of at most two arrows and the target of at most two arrows;
        \item [(G2)] For each arrow $a \in Q_1$, there is at most one arrow $b \in Q_1$ such that $ab \notin I$, and at most one arrow $c \in Q_1$ such that $ca \notin I$;
        \item [(G3)] For each arrow $a \in Q_1$, there is at most one arrow $d \in Q_1$ such that $ad \in I$, and at most one arrow $e \in Q_1$ such $ea \in I$;
        \item [(G4)] $I$ is generated by paths of length 2.
    \end{itemize}
\end{definition}
If we drop Assumption (G3) and replace (G4) by
\begin{itemize}
    \item [(S3)] $I$ is generated by paths of length at least 2,
\end{itemize}
we get the more general notion of \emph{string algebra}.

\begin{examples}\label{gent_a_exs}
\begin{itemize}
    \item [(i)] Let $Q$ be the quiver 
\begin{center}
    $\begin{tikzcd}
  & 2 \arrow[ld,"b"above left] &                         \\
1 &              & 3 \arrow[lu,"a" above right] \arrow[ll,"c" above]
\end{tikzcd}$
\end{center}
 and let $I=\langle ab \rangle$. The algebra $A=kQ/I$ is gentle.
 \item [(ii)] Let $Q$ be the quiver
\begin{center}
    $
\begin{tikzcd}
1 \arrow[rd,"a" left] &              & 3 \arrow[ll,"c" above] \arrow[rd,"d" left] &              & 5 \arrow[ll,"f" above] \\
             & 2 \arrow[ru,"b" right] &                         & 4 \arrow[ru,"e" right] &             
\end{tikzcd}$
\end{center}
and let $I=\langle ab,bc, ca,de,ef,fd \rangle$. The algebra $A=kQ/I$ is gentle.
\end{itemize}    
\end{examples}

Let $A=kQ/I$ be a string algebra. For an arrow $a \in Q_1$, $a^{-1}$ denotes the \emph{formal inverse} of $a$, with $s(a^{-1})=t(a)$ and $t(a^{-1})=s(a)$. The set of formal inverses of the arrows of $Q$ is denoted by $Q_1^{-1}$. The definition of formal inverse is extended to any path $\alpha=a_1 \cdots a_s$ in $Q$, by setting $\alpha^{-1}=a_s^{-1} \cdots a_1^{-1}$.

\begin{definition}[Strings and bands]
 A \emph{string} of length $s\geq 1$ is a walk in the quiver $w=a_1\cdots a_s$ with $a_j \in Q_1 \cup Q_1^{-1}$, and $t(a_j)=s(a_{j+1})$ for which there are no subwalks of the form $aa^{-1}$ or $a^{-1}a$, for some $a \in Q_1$, nor subwalks $ab$ such that $ab \in I$ or $(ab)^{-1} \in I$. For each vertex $i \in Q_0$ there are two \emph{trivial strings} $1_i^{+}$ and $1_i^{-}$ with source and target $x$, for which the formal inverse acts just by swapping the sign. We also consider the zero string $w=0$ of length 0. A string $w=a_1 \cdots a_s$ is \emph{direct} (resp. \emph{inverse}) if $a_i \in Q_1$ (resp. $a_i \in Q_1^{-1}$) for any $i=1,\dots,s$.

For a string $w=a_1\cdots a_s$, we denote by $s(w)$ the source of $a_1$ and by $t(w)$ the target of $a_s$, and we say that $w$ is \emph{cyclic} if $s(w)=t(w)$. A \emph{band} is a cyclic string $b$ such that each power $b^n$ is a string, but $b$ itself is not a proper power of any string. 
 \end{definition}

 \begin{example}\label{ex_band}
     Let $A=kQ/I$ be the gentle algebra given by the quiver
     \begin{center}
         $\begin{tikzcd}
3 \arrow[rr, "d"] \arrow["f"', loop, distance=2em, in=215, out=145] &  & 1 \arrow[rd, "b"'] &                    & 2 \arrow[ll, "a"'] \arrow["e"', loop, distance=2em, in=35, out=325] \\
                                                                    &  &                    & 4 \arrow[ru, "c"'] &                                                                    
\end{tikzcd}$
     \end{center}
     and $I=\langle ab,bc,ca,e^2,f^2\rangle$. Consider $b=a^{-1}ead^{-1}fd$, $v=cead^{-1}$, $w=a^{-1}ea$. Then $b$ is a band, while $v$ and $w$ are strings. In particular, $w$ is an example of a cyclic string which is not a band, as $w^2$ is not a string. 
 \end{example}

\begin{definition}
    Let $v=a_1\dots a_s$ and $w=b_1\dots b_t$ be two strings of length at least 1. The \emph{composition} $vw$ is defined if $vw=a_1\dots a_sb_1\dots b_t$ is a string.
\end{definition}

To define composition with trivial strings and avoid ambiguity in the description of irreducible morphisms between string modules, the following technical definition is required.

\begin{definition}
    Let $A=kQ/I$ be a string algebra. we define two sign functions $\sigma, \epsilon : Q_1 \to \{-1,1\}$ satisfying the following conditions:
\begin{itemize}
    \item [(i)] if $b_1 \neq b_2 \in Q_1$ and $s(b_1)=s(b_2)$, then $\sigma(b_1)=-\sigma(b_2)$;
    \item [(ii)] if $a_1 \neq a_2$ and $t(a_1)=t(a_2)$, then $\epsilon(a_1)=-\epsilon(a_2)$;
    \item [(iii)] if $ab \in I$, then $\sigma(b)=-\epsilon(a)$.
\end{itemize}
The domain of these functions is extended to the set of all strings by setting for any $a \in Q_1$ $\sigma(a^{-1}):=\epsilon(a)$ and $\epsilon(a^{-1}):=\sigma(a)$, for any non-trivial string $w=a_1\dots a_s$, $\sigma(w):=\sigma(a_1)$ and $\epsilon(w):=\epsilon(a_n)$, for any trivial string $1_i$, $\sigma(1_i^{\pm}):=\mp 1$ and $\epsilon(1_i^\pm):=\pm$.
\end{definition} 

\begin{definition}
    Let $v=1_i^\gamma$, $i\in Q_0$, $\gamma \in \{-1,1\}$, and let $w$ be a non-trivial string. The \emph{composition} $vw$ is defined if and only if $s(w)=i$ and $\sigma(w)=-\epsilon(v)=-\gamma$. On the other hand, $wv$ is defined if and only if $t(w)=i$ and $\epsilon(w)=-\sigma(v)=\gamma$.
\end{definition}
\begin{remark}
    If $A$ is gentle, there is always a choice of $\sigma$ and $\epsilon$ such that the concatenation of two strings $v$ and $w$ is defined if and only if $t(v)=s(w)$ and $\sigma(w)=-\epsilon(v)$. However, this is not true for string algebras in general. For an example, see \cite[Remark 1.6]{baur2024geometricmodelmodulecategory}.
\end{remark}
    
\begin{definition}[String modules]
Let $w$ be a string in $A$. The \emph{string module} $M(w)$ associated with $w$ is defined in the following way: the underlying vector space is obtained by replacing each vertex of $w$ by a copy of the field $k$, and the action of an arrow $a \in Q_1$ on $M(w)$ is the identity morphism if $a$ is an arrow of $w$, zero otherwise. If $w=0$, $M(w)$ is the zero module.
\end{definition}
\begin{remark}
\begin{itemize}
    \item [(i)]Let $v,w$ be strings in $A$. Then $M(w)\cong M(v)$ if and only if $v=w$ or $v=w^{-1}$.
    \item [(ii)] If $w$ is a trivial string at vertex $i$ then $M(w)$ is the simple module at $i$.
    \item [(iii)] Indecomposable injective $A$ modules are of the form $I(i) = M(pq^{-1})$, $i \in Q_0$, where $p$ and $q$ are maximal direct paths ending at $i$. The indecomposable projective $A$-modules are of the form $P(i) = M(p^{-1}q)$, $i \in Q_0$, where $p$ and $q$ are maximal direct paths starting from $i$. 
\end{itemize}
    
\end{remark}
\begin{definition}[Band modules]
    Let $b=a_1\cdots a_s$ be a band. The family of \emph{band modules} $M(b,n,\phi)$ associated with $b$, where $n\in \mathbb{N}$ and $\phi$ is an indecomposable automorphism of $k^n$, is defined in the following way: each vertex of $b$ is replaced by a copy of the vector space $k^n$, and the action of an arrow $a \in Q_1$ on $M(b,n,\phi)$ is the identity morphism if $a=a_j$ for $j=1, \dots, n-1$, $\phi$ if $a=a_n$, zero otherwise. The band modules at the mouth of the rank one tubes in the Auslander-Reiten quiver are called \emph{quasi-simple}.  
\end{definition}
\begin{remark}
    Let $b,c$ be bands in $A$. Then $b$ and $c$ define the same (up to isomorphism) family of band modules if and only if $b=c$ up to cyclic permutation and inversion.
\end{remark}
\begin{theorem}[\cite{BR}]
    Let $A$ be a string algebra. All string and band modules are indecomposable. Moreover, every indecomposable $A$-module is either a string module or a band module.
\end{theorem}


Usually strings are drawn as orientations of a Dynkin diagram of type $A$, embedded in the plane so that direct arrows drawn go down and to the right and inverse arrows go down and to the left. One can similarly draw bands as orientations of a type $\Tilde{A}$ Dynkin diagram.

\subsection{Morphisms, extensions and the Auslander-Reiten translation}
\begin{definition}
\begin{itemize}
    \item [(i)] Let $w$ be a string. A \emph{factor string decomposition} of $w$ is a decomposition of the form $w=\Bar{w}_L m_w \Bar{w}_R$, such that $\Bar{w}_L$ is a trivial string or the last letter of $\Bar{w}_L$ lies in $Q_1^{-1}$, and $\Bar{w}_R$ is a trivial string or the first letter of $\Bar{w}_R$ lies in $Q_1$. In this case, $m_w$ is said to be a \emph{factor substring} of $w$. The set of equivalence classes of factor string decompositions of $w$ is denoted by $\mathrm{Fac}(w)$.
    \item [(ii)] Let $v$ be a string. A \emph{substring decomposition} of $v$ is a decomposition of the form $v=\Bar{v}_L m_v \Bar{v}_R$, where $\Bar{v}_L$ is a trivial string or the last letter of $\Bar{v}_L$ lies in $Q_1$, and $\Bar{v}_R$ is a trivial string or the first letter of $\Bar{v}_R$ lies in $Q_1^{-1}$. In this case, $m_v$ is said to be an \emph{image substring} of $v$. The set of equivalence classes of substring decompositions of $v$ is denoted by $\mathrm{Sub}(v)$.
    \item [(iii)] Let $v,w$ be two strings. A pair $(\Bar{w}_L m_w \Bar{w}_R,\Bar{v}_L m_v \Bar{v}_R)$ in $\mathrm{Fac}(w)\times \mathrm{Sub}(v)$ is an \emph{admissible pair} if $m_v=m_w$ or $m_v=m_w^{-1}$. In this case, the string $m_{v,w}:=m_v=m_w$ is called an \emph{overlap between v and w}.
\end{itemize}
    
\end{definition}
\begin{remark}\label{rmk_inj}
   Let $v,w$ be two strings. Each overlap $m_{v,w}$ between $v$ and $w$ defines a morphism $f$ between the string modules $M(w)$ and $M(v)$, given by the composition of the epimorphism $M(w)\to M(m_{v,w})$ and the monomorphism $M(m_{v,w}) \to M(v)$. We will refer to $m_{v,w}$ as the \emph{string associated to} $f$. If $M(m_{v,w})$ is injective, then $\Bar{v}_L=\Bar{v}_R=0$, and so $M(v)=M(m_{v,w})$ is injective. Similarly, if $M(m_{v,w})$ is projective, then $\Bar{w}_L=\Bar{w}_R=0$, and $M(w)=M(m_{v,w})$ is projective.
\end{remark}
\begin{theorem}[\cite{C-B}]
    Let $A$ be a string algebra. Let $v,w$ be two strings and let $M(v)$, $M(w)$ be the corresponding string modules. Then $\mathrm{dim}_k\mathrm{Hom}(M(w),M(v))$ is the number of admissible pairs in $\mathrm{Fac}(w)\times \mathrm{Sub}(v)$.
\end{theorem}

If band modules are involved, the description of morphisms between indecomposable modules over a string algebra is similarly given by overlaps between the corresponding strings or bands. For more details, see \cite{C-B, K}.
    
\begin{definition}
    Let $A=kQ/I$ be a string algebra and let $a \in Q_1$ be an arrow.
    \begin{itemize}
        \item [(i)] Let $p_a$ be the direct string such that $s(p_a)=s(a)$ but it does not start with $a$, and there is no $b \in Q_1$ for which $p_ab$ is a string, i.e. it is \emph{right maximal}. If $s(a)$ has only one outgoing arrow, then $p_a$ is a trivial string at $s(a)$. The \emph{hook of a} is the string $a^{-1}p_a$. 
        \item [(ii)] Let $q_a$ be the direct string such that $t(q_a)=t(a)$ but it does not end with $a$, and there is no $b\in Q_1$ for which $bq_a$ is a string. If $t(a)$ has only one incoming arrow, then $q_a$ is a trivial string at $t(a)$. The \emph{cohook of a} is the string $aq_a^{-1}$.  
    \end{itemize}
    Let $w$ be a string. The string $f_l(w)$ is defined as follows:
     \begin{itemize}
        \item if there exists $a \in Q_1$ such that $aw$ is a string, then $f_l(w):=p_a^{-1}aw$, i.e. $f_l(w)$ is obtained by adding the hook of $a$ on the left of $w$;
        \item if there is no $a \in Q_1$ such that $aw$ is a string and $w$ is not direct, then $w$ must be of the form $q_aa^{-1}\Bar{w}_R$, for some $a \in Q_1$, and $f_l(w):=\Bar{w}_R$, i.e. $f_l(w)$ is obtained by removing the cohook of $a$ on the left of $w$;
        \item if there is no $a \in Q_1$ such that $aw$ is a string and $w$ is direct, then $f_l(w):=0$.
    \end{itemize}
    Similarly, the string $f_r(w)$ is defined as follows:
    \begin{itemize}
        \item if there exists $a \in Q_1$ such that $wa^{-1}$ is a string, then $f_r(w):=wa^{-1}p_a$, i.e. $f_r(w)$ is obtained by adding the hook of $a$ on the right of $w$;
        \item if there is no $a \in Q_1$ such that $wa^{-1}$ is a string and $w$ is not inverse, then $w$ must be of the form $\Bar{w}_Laq_a^{-1}$, for some $a \in Q_1$, and $f_r(w):=\Bar{w}_L$, i.e. $f_r(w)$ is obtained by removing the cohook of $a$ on the right of $w$; 
        \item if there is no $a \in Q_1$ such that $wa^{-1}$ is a string and $w$ is inverse, then $f_r(w):=0$.
    \end{itemize}
    In the same way, $g_l(w)$ (resp. $g_r(w)$) is defined as the string obtained from $w$ by either adding a cohook or removing a hook on the left (resp. on the right). 
\end{definition}


\begin{theorem}[\cite{BR}]\label{AR_seq_string_alg}
    Let $A$ be a string algebra and $w$ a string over $A$. 
    \begin{itemize}
        \item [(i)] There are at most two irreducible morphisms starting at $M(w)$, with possible targets $M(f_l(w))$ and $M(f_r(w))$.
        \item [(ii)] If $M(w)$ is not an injective module, assuming without loss of generality that $f_l(w)\neq 0$, then the Auslander-Reiten sequence starting at $M(w)$ is of the form
        \begin{equation*}
            0 \to M(w) \to M(f_l(w)) \oplus M(f_r(w)) \to M(f_r(f_l(w))) \to 0.
        \end{equation*}
        In particular, $\tau^{-1}(M(w))=M(f_r(f_l(w)))$.
        \item [(iii)] If $M(w)$ is not a projective module, assuming without loss of generality that $g_l(w)\neq 0$, then the Auslander-Reiten sequence ending at $M(w)$ is of the form
        \begin{equation*}
            0 \to  M(g_r(g_l(w))) \to M(g_l(w)) \oplus M(g_r(w))  \to M(w) \to 0.
        \end{equation*}
        In particular, $\tau(M(w))=M(g_r(g_l(w)))$.
    \end{itemize}
\end{theorem}

From now on, the term 'band module' will refer to a quasi-simple band module.

\begin{figure}
    \centering
    \includegraphics[width=1\linewidth]{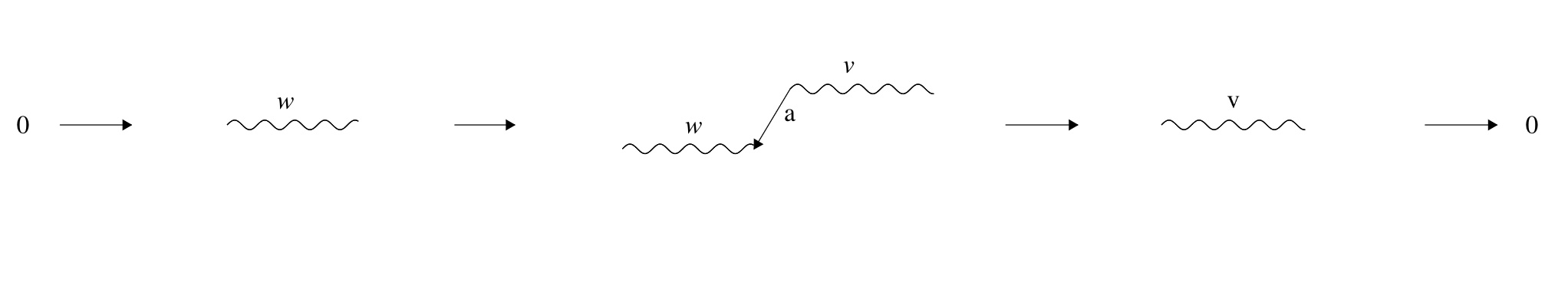}
    \caption{An arrow extension between $M(w)$ and $M(v)$.}
    \label{arrow_ext}
\end{figure}

\begin{figure}
    \centering
    \includegraphics[width=1\linewidth]{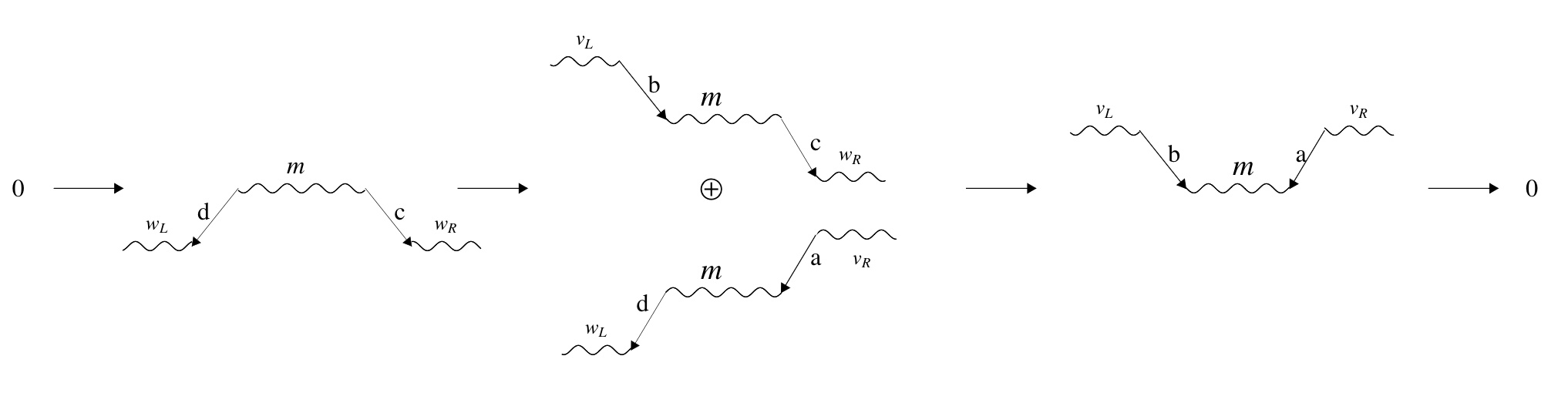}
    \caption{An overlap extension between $M(w)$ and $M(v)$ with overlap $m$.}
    \label{overlap_ext}
\end{figure}

\begin{definition}[{\cite[Definition 3.1]{CPS}}]\label{def_ext_gent_alg}
  Let $v$ and $w$ be string. 
  \begin{itemize}
      \item [(1)] \textbf{(Arrow extension)} If there exists $a \in Q_1$ such that $u = wa^{-1}v$ is a string then there is
a non-split short exact sequence (Figure \ref{arrow_ext})
\begin{center}
    $0 \to M(w) \to M(u) \to M(v) \to 0$.
\end{center}
\item [(2)]\textbf{(Overlap extension)} Suppose that $v = v_Lb ma^{-1}v_R$ and $w = w_Ld^{-1}mcw_R$ with $a, b, c, d \in Q_1$ and $m, v_L, v_R, w_L, w_R$ (possibly trivial) strings such that
\begin{itemize}
    \item [(i)] if $a=\emptyset$, then $c \neq \emptyset$;
    \item [(ii)] if $b=\emptyset$, then $d \neq \emptyset$; 
    \item [(iii)] if $m$ is a trivial string, then $ca \in I$ and $bd \in I$ (whenever they exist, subject to the constraints above).
 \end{itemize}
Then there is a non-split short exact sequence
\begin{center}
    $0 \to M(w) \to M(u) \oplus M(u') \to M(v) \to 0$,
\end{center}
where $u = v_Lbmc w_R$ and $u' = w_Ld^{-1}ma^{-1}v_R$ (Figure \ref{overlap_ext}).
  \end{itemize}
\end{definition}

\begin{remark}\label{rmk_rel_ext}
  If $A$ is a gentle algebra and $v,w$ are strings as in \ref{def_ext_gent_alg} (2), then $ca \in I$ and $bd \in I$.  
\end{remark}

\begin{theorem}[{\cite[Theorem A]{CPS}}]\label{thm_A}
    Let $A$ be a gentle algebra, and $v$ and $w$ strings. The collection of arrow and overlap extensions between $M(w)$ and $M(v)$ form a basis of $\mathrm{Ext}_A^1(M(v),M(w))$.
\end{theorem}

\begin{example}
  Let $A$ be the algebra of Example \ref{gent_a_exs} (ii). Let $w=c^{-1}d$ and $v=bf^{-1}$. Then, by Theorem \ref{thm_A}, the non-split short exact sequence
  \begin{center}
   $0\to M(c^{-1}d)\to M(bd)\oplus M(c^{-1}f^{-1})\to M(bf^{-1})\to 0$  
  \end{center}
  generates $\mathrm{Ext}_A^1(M(v),M(w))$. 
\end{example}

If band modules are involved, there are no arrow extensions, only overlap extensions. For more details, see \cite[Theorems B-E]{CPS}.

\begin{remark}
    If $A$ is a non-gentle string algebra and $v, w$ are strings or bands, then overlaps between $v$ and $w$ no longer provide a basis for $\mathrm{Ext}_A^1(M(v), M(w))$, as noted in \cite[Remark 2.7]{Z}. For instance, consider the algebra $A = kQ/I$, where $Q = 1 \xrightarrow{a} 2 \xrightarrow{b} 3 \xrightarrow{c} 4$ and $I = \langle abc \rangle$. $A$ is a string algebra that is not gentle. Take the two strings $v = ab$ and $w = bc$. There is an overlap between $v$ and $w$ given by $b$. However, $\mathrm{Ext}_A^1(M(v), M(w)) = 0$ since $M(v)$ is projective.
\end{remark}

\begin{definition}\label{rigid_def}
    A module $L$ over a finite-dimensional algebra is called \emph{rigid} if $\mathrm{Ext}^1(L,L)=0$.
\end{definition}

\begin{example}\label{nonrigid_band}
    A band module is never rigid, as the Auslander-Reiten translation acts on it as the identity.
\end{example}

\subsection{Gentle algebras arising from triangulations of surfaces}
Let $(S,M)$ be an unpunctured bordered surface with marked points. That is, $S$ is an oriented surface with boundary $\delta S$, and $M \subset \delta S$ is a non-empty finite set of marked points on $\delta S$ that intersects each connected component of the boundary $\delta S$. An \emph{arc} in $(S,M)$ is a homotopy class of curves connecting marked points. An ideal triangulation of $(S,M)$ is a maximal collection of non-crossing arcs.

Let $T$ be an ideal triangulation of $(S,M)$. Naturally associated with $T$ is an algebra $A(T)$ given as the Jacobian algebra of a quiver with potential arising from the triangulation \cite{LF}. By \cite{ABCP}, $A(T)$ is a gentle algebra. Moreover, arcs in $(S,M)$ are in bijection with string modules over $A(T)$, while band modules are in bijection with homotopy classes of closed curves.
\begin{definition}
    $A(T)$ is called the \emph{gentle algebra of the triangulation $T$}.
\end{definition}
\begin{remark}\label{rmk_rel_jac}
    Let $A(T)=kQ/I$ be the gentle algebra of $T$. Then any oriented cycle in $Q$ is a 3-cycle with full relations, that is, the composition of any two arrows in it is in $I$.
\end{remark}

\begin{example}\label{ex_alg_triang}
    The algebra $A$ in Example \ref{gent_a_exs} (ii) is the gentle algebra associated with the triangulation $T'$ of the octagon in Figure \ref{fig:ex_typeB}.
\end{example}

\section{Cluster character, $F$-polynomial and $\bold{g}$-vector of a module}\label{section_cc}
 Grassmannians of vector spaces are projective varieties whose points parametrize subvector spaces of a given dimension. Quiver Granssmannians generalize this notion: given a finite-dimensional algebra $A$, they are projective varieties whose points parametrize submodules of a given dimension vector, that is, the graded dimension of the underlying vector space.

\begin{definition}\label{defi::grassmannian}
Let $A$ be a finite-dimensional algebra. Let $\bold{e} \in \mathbb{Z}_{\geq 0}^{Q_0}$ and $L$ be an $A$-module. The \emph{quiver Grassmannian of $L$ with dimension vector $\bold{e}$} is the set $\mathrm{Gr}_{\bold{e}}^A(L)$ of all submodules of $L$ of dimension vector $\bold{e}$. If the algebra is clear from the context, we simply denote it as $\mathrm{Gr}_{\bold{e}}(L)$.
\end{definition}

\begin{examples}\label{exam::Grassmannians}
 \begin{itemize}
  \item[(i)] If the quiver $Q$ has only one vertex and no arrows, then the $kQ$-modules are just vector spaces, and their quiver Grassmannians are just usual Grassmannians.
  \item [(ii)] Let $A$ be the algebra in Example \ref{gent_a_exs} (i). Let $w=a^{-1}c b^{-1}$ and let $L=M(w)$ be the corresponding string module. The following table lists the dimension vectors for which the quiver Grassmannian of $L$ is non-empty, and for each of these vectors gives a variety isomorphic to the corresponding quiver Grassmannian.
    \[
        \begin{tabular}{|c||c|c|c|c|c|c|c|} 
          \hline
          $\bold{e}$   & $(0,0,0)$ & $(1,0,0)$ & $(0,1,0)$ & $(1,1,0)$ & $(1,2,0)$ & $(1,1,1)$ & $(1,2,1)$ \\
          \hline 
          $\mathrm{Gr}_{\bold{e}}(L)$ & point & point & point & $\mathbb{P}^1$ & point & point & point \\
          \hline
        \end{tabular}
    \] 
  \end{itemize}
  \end{examples}
  
\begin{definition}\label{defi::F-polynomial}
 Let $L$ be an $A$-module. Its \emph{$F$-polynomial} is 
 $$ F_L(\bold{y}) := \displaystyle\sum_{\bold{e}\in\mathbb{Z}_{\geq 0}^{Q_0}} \chi\big( \mathrm{Gr}_{\bold{e}}(L) \big) \bold{y}^{\bold{e}},
 $$
 where $\bold{y}^{\bold{e}} = \displaystyle\prod_{i\in Q_0} y_i^{e_i}$ and $\chi$ is the Euler-Poincar\'e characteristic.
\end{definition}
\begin{remark}\label{rmk_EC}
Let $Q$ be a quiver. A \emph{successor closed subquiver} of $Q$ is a subquiver $T$ such that if $i \in T_0$ is a vertex of $T$ and $a:i\to j \in Q_1$ is an arrow of $Q$, then $a \in T_1$ is an arrow of $T$. 
  If $A$ is a string algebra and $L$ is a string or a band module, then $\chi(\mathrm{Gr}_{\bold{e}}(L))$ is the number of successor closed subquivers of the diagram of type $A$ or type $\Tilde{A}$ representing $L$. For more details, see \cite{CI,JS, Haupt}. The \emph{dimension} of a successor closed subquiver $T$, denoted by $\bold{dim}T$, is the non-negative integer vector whose $i$-th coordinate is the number of times that $i \in Q_0$ appears in $T$.      
\end{remark}
\begin{examples}\label{ex_F_pol}
\begin{itemize}
    \item [(i)] Consider the module $L$ of Example \ref{exam::Grassmannians} (ii). The module $L$ can be represented with the following diagram:
    \begin{center}
        \begin{tikzcd}
  & 3 \arrow[ld, "a"'] \arrow[rd, "c"] &   & 2 \arrow[ld, "b"'] \\
2 &                                              & 1 &                       
\end{tikzcd}
    \end{center}
    It follows that $F_L=1+y_1+y_2+2y_1y_2+y_1y_2^2+y_1y_2y_3+y_1y_2^2y_3$.
    \item [(ii)] Consider the algebra $A$ of Example \ref{ex_band}. Let $L$ be the (quasi-simple) band module corresponding to the band $b=a^{-1}ead^{-1}fd$. It can be represented with the following diagram:
    \begin{center}
        $\begin{tikzcd}
  & 2 \arrow[ld, "a"'] \arrow[rd, "e"] &                   &   &                                    &                                     \\
1 &                                    & 2 \arrow[rd, "a"] &   & 3 \arrow[ld, "d"'] \arrow[rd, "f"] &                                     \\
  &                                    &                   & 1 &                                    & 3 \arrow[lllllu, "d", bend left=49]
\end{tikzcd}$
    \end{center}
    One can compute $F_L=1+2y_1+y_1^2+y_1y_2+y_1y_3+y_1^2y_2+y_1^2y_3+y_1^2y_2^2+y_1^2y_3^2+y_1^2y_2y_3+y_1^2y_2y_3^2+y_1^2y_2^2y_3+y_1^2y_2^2y_3^2$.
\end{itemize}
    
\end{examples}

\begin{definition}\label{def_g_vector}
    Let $L$ be an $A$-module. Let
    \begin{center}
        $0 \to L \xrightarrow[]{} I_0 \xrightarrow[]{} I_1$
    \end{center}
    be a minimal injective presentation of $L$ in $\mathrm{mod}\hspace{0.05cm}A$, where $I_0=\displaystyle\bigoplus_{i\in Q_0}I(i)^{a_i}$ and $I_1=\displaystyle\bigoplus_{i\in Q_0}I(i)^{b_i}$. Then the $\bold{g}$-vector of $L$ is the vector $\bold{g}_L \in \mathbb{Z}^{Q_0}$ whose $i$-th coordinate is given by
    \begin{center}
        $(\bold{g}_L)_i:=b_i-a_i$.
    \end{center}
\end{definition}

 For string modules over gentle algebras, $\bold{g}$-vectors are easy to compute using string combinatorics. 

\begin{theorem}[{\cite[Corollary 1.44]{PPP}}]\label{thm_g_vect}
    Let $A=kQ/I$ be a gentle algebra and let $w$ be a string. Then 
    \begin{align*}
        \bold{g}_{M(w)}=\bold{b}+\bold{r}-\bold{a},
    \end{align*}
    where
     \begin{itemize}
  \item $\bold{a}=(a_i)_{i\in Q_0}$ where $a_i$ is the number of times the vertex $i$ is found at the bottom of $w$;
  \item $\bold{b}=(b_i)_{i \in Q_0}$ where $b_i$ is the number of times the vertex $i$ is found at the top of $w$, but not at its start or end;
        \medskip
  \item 
  		\(
         \bold{r} = \begin{cases}
                  \bold{e}_{s(a)}+\bold{e}_{s(b)}  & \textrm{if $awb^{-1}$ is a string, $a,b \in Q_1$;}\\
                   \bold{e}_{s(a)}   & \textrm{if $aw$ is a string, $a \in Q_1$, and there is no $b \in Q_1$ such that $wb^{-1}$ is a string;} \\
                   \bold{e}_{s(b)}   & \textrm{if $wb^{-1}$ is a string, $b\in Q_1$, and there is no $a \in Q_1$ such that $aw$ is a string;} \\
                   0 & \textrm{if there are no $a,b \in Q_1$ such that $aw$ or $wb^{-1}$ is a string,}
                  \end{cases}
        \) \\
        \medskip
        where $(\bold{e}_i)_{i \in Q_0}$ is the standard basis of~$\mathbb{Z}^{Q_0}$.
 \end{itemize}
\end{theorem}

\begin{examples}
    \begin{itemize}
        \item [(i)] Let $A$ be the algebra in Example \ref{gent_a_exs} (i) and let $L=M(1_1)$ be the simple in 1.
        \begin{center}
            $\begin{tikzcd}
2 \arrow[rd, "b"', dotted] &   & 3 \arrow[ld, "c", dotted] \\
                           & 1 &                          
\end{tikzcd}$
        \end{center}
        By Theorem \ref{thm_g_vect}, $\bold{g}_L=\begin{pmatrix}
            -1\\1\\1
        \end{pmatrix}$.
    \end{itemize}
    \item [(ii)] Let $A$ be the algebra in Example \ref{gent_a_exs} (ii) and let $L=M(c^{-1}d)$.
    \begin{center}
        \begin{tikzcd}
\times \arrow[rd, dotted] &   & 3 \arrow[ld,"c"'] \arrow[rd,"d"] &   & \times \arrow[ld, dotted] \\
                          & 1 &                         & 4 &                          
\end{tikzcd}
    \end{center}
     By Theorem \ref{thm_g_vect}, $\bold{g}_L=\begin{pmatrix}
            -1\\0\\1\\-1 \\0
        \end{pmatrix}$.
\end{examples}

To a quiver $Q$ is associated a skew-symmetric matrix $B(Q)=(b_{ij})$, where $b_{ij}=|\{a:j\to i \in Q_1\}|-|\{a:i\to j \in Q_1\}|$, $i\neq j$, and $b_{ii}=0$. $B(Q)$ does not detect the presence of loops in $Q$.

\begin{definition}\label{cc_map}
    Let $A=kQ/I$ be a finite-dimensional algebra. Let $B=B(Q)$, and let $n=|Q_0|$. Let $L$ be an $A$-module. The \emph{cluster character of $L$}, also known as \emph{Caldero-Chapoton map}, is the Laurent polynomial
    \begin{center}
        $CC(L)=\displaystyle\sum_{\bold{e}\in \mathbb{Z}_{\geq 0}^n}\chi(\mathrm{Gr}_{\bold{e}}(L))\bold{x}^{\,B\bold{e}+\bold{g}_L}\,\bold{y}^{\bold{e}}\in \mathbb{Z}[y_1,\dots,y_n][x_1^{\pm 1},\dots,x_n^{\pm 1}]$.
    \end{center}
\end{definition}

\begin{remark}
    One can obtain $F_L$ by specializing $CC(L)$ at $x_1=\cdots=x_n=1$. On the other hand, $\bold{g}$ is the multi-degree of $CC(L)$ with respect to the $\mathbb{Z}^n$-grading in $\mathbb{Z}[y_1,\dots,y_n][x_1^{\pm 1},\dots,x_n^{\pm 1}]$ given by $\mathrm{deg}(x_i)=\bold{e_i}$ and $\mathrm{deg}(y_j)=-\bold{b}_j$, where $\bold{e_i}$ is the $i$-th vector of the standard basis of $\mathbb{Z}^n$, and $\bold{b}_j$ is the $j$-th column of $B$.
\end{remark}

\section{Cluster multiplication formula for acyclic quivers}\label{s_acyc}

Let $A$ be an algebra and $L,N$ modules over $A$. We use the following notation:
\begin{center}
    $[L,N]:=\text{dim}_k\mathrm{Hom}_A(L,N)$ \hspace{1cm} and \hspace{1cm} $[L,N]^1:=\text{dim}_k\mathrm{Ext}_A^1(L,N)$.
\end{center}
\begin{definition}

An element  $\xi : 0 \to X \to Y \to S \to 0 \in\mathrm{Ext}^1(S,X)$ is  \emph{generating} if $\mathrm{Ext}^1(S,X)=\mathbb{C}\xi$. 
\end{definition}
In other words  $\xi\in\mathrm{Ext}^1(S,X)$ is generating if either $[S,X]^1=0$ and $\xi=0$, or $[S,X]^1=1$ and $\xi\neq 0$.  

 \vspace{0.2cm}

Let $Q$ be an acyclic quiver. Let $X$ and $S$ be $kQ$-modules such that $[S, X]^1 = 1$, and let
\begin{center}
    $\xi:0 \to X \to Y \to S \to 0$
\end{center} be a generating extension. In \cite[Lemma 27]{CEHR}, it is proved that the following submodules of $X$ and $S$ are well defined:
$$
\begin{array}{cc}
\overline{X}:=\textrm{max}\{N\subset X|\, [S,X/N]^1=1\} \subset X,& \underline{S}:=\textrm{min}\{N\subset S|\,[N,X]^1=1\}\subseteq S.
\end{array}
$$
In other words, the submodule $\overline{X}$ is the maximal submodule of $X$ such that the ``push-out'' sequence

\begin{center}
    $\begin{tikzcd}
        \xi:0 \arrow[r] & X \arrow[r] \arrow[d,two heads, "p"] & Y \arrow[r] \arrow[d,two heads] & S \arrow[r] \arrow[d,"=" right] & 0
        \\ p_\ast\xi : 0 \arrow[r] & X/\overline{X} \arrow[r] & \Bar{Y} \arrow[r] & S \arrow[r] & 0
    \end{tikzcd}$
\end{center}

does not split. Dually, the submodule $\underline{S}\subseteq S$ is the minimal submodule such that the ``pull-back'' sequence

\begin{center}
    $\begin{tikzcd}
        i^\ast \xi:0 \arrow[r] & X \arrow[r] \arrow[d, "="] & \Tilde{Y} \arrow[r] \arrow[d,hook] & \underline{S} \arrow[r] \arrow[d,hook,"i"] & 0
        \\ \xi : 0 \arrow[r] & X \arrow[r] & \Bar{Y} \arrow[r] & S \arrow[r] & 0
    \end{tikzcd}$
\end{center}

does not split. If $\xi$ is almost split, then this description implies that $\underline{S}=S$ and $\overline{X}=0$.
\begin{example}\label{ex_xs_sx}
    Let $Q$ be the quiver $1 \gets 2 \gets 3 \gets 4 \gets 5$. 
    We consider the short exact sequence 
    \begin{center}
        $\xi: 0 \to X=\begin{smallmatrix}
                 3\\2\\1
             \end{smallmatrix} \to Y=\begin{smallmatrix}
                 5\\4\\ 3\\2\\1
             \end{smallmatrix} \oplus \begin{smallmatrix}
                3
             \end{smallmatrix}\to S=\begin{smallmatrix}
                 5\\4\\ 3
             \end{smallmatrix}\to 0$.
    \end{center}
    Then $\overline{X}=\begin{smallmatrix}
        1
    \end{smallmatrix}$ and $\underline{S}=\begin{smallmatrix}
        4\\3
    \end{smallmatrix}$.
\end{example}

\vspace{0.2cm}

\begin{definition}
  A generating extension $\xi : 0 \to X \to Y \to S \to 0$ is called a \emph{generalized almost split} sequence if $\underline{S}=S$ and $\overline{X}=0$.  
\end{definition}

\begin{remark}
    An almost split sequence is generalized almost split.
\end{remark}

\begin{example}\label{counterexample}
    We exhibit an example of generalized almost split sequence which is not almost split for a type $A$ quiver. This is a counterexample for \cite[Example 30]{CEHR} which asserts that, for quivers of type $A$, generalized almost split sequences are almost split. 

    Let $Q$ be the quiver
    \begin{center}
        $\begin{tikzcd}
& 2 \arrow[dl] \arrow[d] \\
1  & 3 \arrow[d] & 5 \arrow[dl]  \\
&  4 
\end{tikzcd}$
    \end{center}
The Auslander-Reiten quiver of $KQ$ is:
    \begin{center}
        $\begin{tikzcd}
                        & \begin{smallmatrix}
                            5\\4
                        \end{smallmatrix} \arrow[rd]            &                            & \begin{smallmatrix}
                            3
                        \end{smallmatrix} \arrow[rd]                &                            & \begin{smallmatrix}
                            2\\1
                        \end{smallmatrix} \arrow[rd] &   \\
\begin{smallmatrix}
    4
\end{smallmatrix} \arrow[ru] \arrow[rd] &                          & \begin{smallmatrix}
    35\\4
\end{smallmatrix} \arrow[ru] \arrow[rd]  &                             & \begin{smallmatrix}
    2\\13
\end{smallmatrix} \arrow[ru] \arrow[rd]  &               & \begin{smallmatrix}
    2
\end{smallmatrix} \\
                        & \begin{smallmatrix}
                            3\\4
                        \end{smallmatrix} \arrow[ru] \arrow[rd] &                            & \begin{smallmatrix}
                            2\\135\\4
                        \end{smallmatrix} \arrow[ru] \arrow[rd] &                            & \begin{smallmatrix}
                            2\\3
                        \end{smallmatrix} \arrow[ru] &   \\
                        &                          & \begin{smallmatrix}
                            2\\13\\4
                        \end{smallmatrix} \arrow[ru] \arrow[rd] &                             & \begin{smallmatrix}
                            2\\35\\4
                        \end{smallmatrix} \arrow[ru] \arrow[rd] &               &   \\
                        & \begin{smallmatrix}
                            1
                        \end{smallmatrix} \arrow[ru]             &                            & \begin{smallmatrix}
                            2\\3\\4
                        \end{smallmatrix} \arrow[ru]              &                            & \begin{smallmatrix}
                            5
                        \end{smallmatrix}             &  
\end{tikzcd}$
    \end{center}
We consider the short exact sequence
\begin{center}
    $\xi:0\to \begin{smallmatrix}
        3\\4
    \end{smallmatrix}\to \begin{smallmatrix}
        3 
    \end{smallmatrix}\oplus\begin{smallmatrix}
        2\\3\\4
    \end{smallmatrix}\to \begin{smallmatrix}
        2\\3
    \end{smallmatrix}\to 0$.
\end{center}
We have that $\xi$ is generalized almost split, but not almost split.
\end{example}

\begin{theorem}[\cite{CEHR}]\label{thm_start}
Let $X$, $S$ be $kQ$-modules such that $[S,X]^1=1$. Let $\xi: 0 \to X \to Y \to S \to 0 \in \mathrm{Ext}^1(S,X)$ be a generating extension with a middle term $Y$. Then for any $\bold{e}\in \mathbb{Z}_{\geq0}^{Q_0}$
    \begin{equation}\label{thm_start_eq}
        \chi(Gr_{\bold{e}}(X \oplus S))=\chi(Gr_{\bold{e}}(Y))\chi( Gr_{\bold{e}-\textbf{dim}\underline{S}}(\overline{X} \oplus S/\underline{S})).
    \end{equation}
\end{theorem}

Let X, S be $kQ$-modules such that $[S,X]^1=1$. Then, there exists an exact sequence $0 \to X/\overline{X} \to \tau \underline{S} \to I \to 0$,
where $I$ is either injective or zero \cite[Lemma 31]{CEHR}. Let $|Q_0|=n$, let $I = I(1)^{f_1} \oplus I(2)^{f_2} \oplus \cdots \oplus I(n)^{f_n}$  be the indecomposable decomposition of $I$, and let $\bold{f} = (f_1,\cdots,f_n)$. Theorem \ref{thm_start} leads to the following multiplication formula for cluster characters:

\begin{theorem}[{\cite[Theorem 67]{CEHR}}]\label{thm67}
  Let $\xi \in \mathrm{Ext}^1(S,X)$ be a generating extension with middle term $Y$. Then
    \begin{equation}\label{eq_thm67}
        CC(X)CC(S) = CC(Y) + \bold{y}^{\textbf{dim} \underline{S}} CC(\overline{X} \oplus S/\underline{S} )\bold{x}^{\bold{f}},
    \end{equation}
  where $\bold{x}=(x_1,\dots,x_n)$ and $\bold{y}=(y_1,\dots,y_n)$. Moreover, if $\mathrm{Ext}^1(X,S) = 0$, and both $X$ and $S$ are rigid and indecomposable, then the formula \ref{eq_thm67} is an exchange relation between the cluster variables $CC(X)$ and $CC(S)$ for the cluster algebra with principal coefficients in the initial seed whose exchange matrix is $B(Q)$.
 
\end{theorem}
\begin{example}
    Let $Q, \xi, X, S$ as in Example \ref{ex_xs_sx}. By Theorem \ref{thm67}, we have
    \begin{equation}\label{eq_exchange}
        CC(X)CC(S)=CC(\begin{smallmatrix}
            3\\2\\1
        \end{smallmatrix})CC(\begin{smallmatrix}
            5\\4\\3
        \end{smallmatrix})=CC(\begin{smallmatrix}
            5\\4\\3\\2\\1
        \end{smallmatrix}\oplus\begin{smallmatrix}
            3
        \end{smallmatrix})+y_3y_4CC(\begin{smallmatrix}
            1
        \end{smallmatrix}\oplus\begin{smallmatrix}
            5
        \end{smallmatrix})
    \end{equation}
     
    In this case $\bold{f}$ is zero in \ref{eq_thm67}, since $\tau \underline{S}=X/\overline{X}$. Since both $X$ and $S$ are rigid and indecomposable, \ref{eq_exchange} is an exchange relation for the cluster algebra $\mathcal{A}_\bullet(B)$ of type $A_5$, where $B=B(Q)$. 
\end{example}

\section{The multiplication formula for gentle algebras}\label{s_gent_alg_cc}

The definitions of $\overline{X}$ and $\underline{S}$, as well as the proof of Theorem \ref{thm_start}, heavily rely on the fact that the algebra $kQ$ is hereditary, meaning $\mathrm{Ext}_{kQ}^i(-,-) = 0$ for any $i > 1$. To generalize Theorem \ref{thm_start} in the new setting, the idea is to define $\overline{X}$ and $\underline{S}$ differently, as follows:
\begin{definition}
Let $A=kQ/I$ be a gentle algebra.  Let $X$, $S$ be $A$-modules such that $[S,X]^1=1$.
   \begin{center} $\overline{X} := \text{ker}(f) \subset X$; \hspace{2cm}$\underline{S} := \text{im}(g) \subseteq S$ \end{center}
where $f: X \to \tau S$ is the non-zero morphism from $X$ to $\tau S$ that does not factor through an injective $A$-module, and $g: \tau^{-1} X \to S$ is the non-zero morphism from $\tau^{-1} X$ to $S$ that does not factor through a projective $A$-module. 
\end{definition}
\begin{remark}
  By the Auslander-Reiten formulas (cf. \cite{libro_blu}, IV.2.13), $\overline{X}$ and $\underline{S}$ are well-defined $A$-modules. Moreover, according to \cite[Lemma 31]{CEHR}, the new definitions coincide with the original ones when $A$ is hereditary, that is, $Q$ has no oriented cycles and $I=0$ (cf. \cite{libro_blu}, VII.1.7).  
\end{remark}

With these new definitions, we obtain the first main result in this work. To state the theorem, we need to recall the $\mathrm{Ext}$-order between modules over a finite-dimensional algebra.
 
\begin{definition}
Let $A$ be a finite-dimensional algebra, and let $M,N$ be $A$-modules. We say that $M \leq_{\mathrm{Ext}} N$ if there exist $A$-modules $M_1,\dots, M_k$ such that for every $i$ there exists a short exact sequence
\begin{center}
    $0 \to U_i \to M_{i-1} \to V_i \to 0$,
\end{center}
such that $M_1=M$, $M_k=N$, $M_i \cong U_i \oplus V_i$.
\end{definition}

\begin{theorem}\label{dec_qg_gent_alg}
     Let $A=kQ/I$ be a gentle algebra.  Let $X$, $S$ be $A$-modules such that $[S,X]^1=1$. Let $\xi: 0 \to X \to Y \to S \to 0 \in \mathrm{Ext}^1(S,X)$ be a generating extension with a middle term $Y$. Then for any $\bold{e}\in \mathbb{Z}_{\geq0}^{Q_0}$
    \begin{equation}\label{dec_qg_gent_alg_eq}
        \chi(Gr_{\bold{e}}^A(X \oplus S))=\chi(Gr_{\bold{e}}^A(Y))+\chi( Gr_{\bold{e}-\textbf{dim}\underline{S}}^{B}(M)),
    \end{equation} 
    where $M$ is the $\leq_{\mathrm{Ext}}$-minimum extension between $S/\underline{S}$ and $\overline{X}$ in a finite-dimensional algebra $B=kQ'/I\supseteq A$. In particular, if $A$ is the gentle algebra of a triangulation of an unpunctured surface with marked points, then $B=A$. Moreover, $M <_{\mathrm{Ext}}\overline{X}\oplus S/\underline{S}$ if and only if $\mathrm{Ext}_A^1(\underline{S},X/\overline{X})=0$.
\end{theorem}

\begin{remark}\label{M_rigid}
    The module $M$ is rigid.
\end{remark}

\begin{remark}
    The algebra $B$ explicitly constructed in the proof of Theorem \ref{dec_qg_gent_alg} is generally no longer gentle. See Example \ref{ex_dec_thm} (iii). 
\end{remark}

\begin{remark}
    If $A$ is hereditary, by \cite[Lemma 31]{CEHR} $\mathrm{Ext}^1(\underline{S},X/\overline{X}) \neq 0$ for any $X,S$ such that $[S,X]^1=1$, so we recover the statement of Theorem \ref{thm_start}.
\end{remark}

\begin{proof}
By the additivity of $\mathrm{Ext}^1$, we can consider without loss of generality the case where $X$ and $S$ are indecomposable. In fact, if $X=\displaystyle\bigoplus_{i=1}^s X_i$ and $S=\displaystyle\bigoplus_{j=1}^t S_j$, then $1=[S,X]^1=\displaystyle\sum_{i,j}[S_j,X_i]^1$ implies $[S_h,X_k]^1=1$ for some $1\leq h \leq t$, $1 \leq k \leq s$ and $[S_j,X_i]^1=0$ for any $i\neq k,j\neq h$, so $\overline{X}=\overline{X_k}\displaystyle\bigoplus_{i\neq k}X_i$, $\underline{S}=\underline{S_h}$, and $S/\underline{S}=S_h/\underline{S_h} \displaystyle\bigoplus_{j\neq h}S_j$.

We show the proof for $X,S$ both string modules. If one or both are bands, the argument is the same with minor adaptations. Let $X=M(w)$ and $S=M(v)$ be string modules such that $[S,X]^1=1$, and let $\xi: 0 \to X \to Y \to S \to 0 \in \mathrm{Ext}^1(S,X)$ be a generating extension. 

We first consider the case where $\xi$ is an overlap extension. Since $[S,X]^1=1$, by Theorem \ref{thm_A} there is only one overlap between $v$ and $w$. Consequently, there is only one overlap between the string $v'=g_r(g_l(v))$ corresponding to $\tau S$ and $w$, which provides a map between the corresponding modules that does not factor through an injective module. Furthermore, there is only one maximum length overlap $m_\mathrm{max}$ in $X$ between $v'$ and $w$. In fact, every other overlap must be a substring of $m_\mathrm{max}$ in $X$; otherwise, it would also be an overlap between $v$ and $w$, contradicting the hypothesis that $[S,X]^1=1$.
We claim that the morphism $f$ associated with $m_\mathrm{max}$ does not factor through an injective module. In fact, let $g$ be a morphism that factors through an injective module $I=M(u)$. Then there are overlaps $m_{u,w}$ between $u$ and $w$, and $m_{v',u}$ between $v'$ and $u$, and $g$ is the morphism associated with the overlap $m_{u,w} \cap m_{v',u}$ between $v'$ and $w$, which is a proper substring of $m_\mathrm{max}$. Otherwise, since $I$ is injective and $m_\mathrm{max}$ is maximal, it would follow that $u=m_{u,w} \cap m_{v',u}=m_\mathrm{max}$. Consequently, $M(m_\mathrm{max})$ would be injective, and by Remark \ref{rmk_inj}, $\tau S$ would also be injective, which leads to a contradiction. Hence, $g \neq f$.
Similarly, the morphism from $\tau^{-1}X$ to $S$ associated with the maximum overlap in $S$ between the corresponding strings does not factor through a projective module.


It follows that $\overline{X}=M(\alpha_L)\oplus M(\alpha_R)$ and $S/\underline{S}=M(\beta_L)\oplus M(\beta_R)$, where $\alpha_L$ and $\alpha_R$ are zero strings or image substrings of $w_L$ and $w_R$ such that $s(w_L)=s(\alpha_L)$ and $t(w_R)=t(\alpha_R)$, and $\beta_L$ and $\beta_R$ are zero strings or factor substrings of $v_L$ and $v_R$ such that $s(v_L)=s(\beta_L)$ and $t(v_R)= t(\beta_R)$. See Figure \ref{non-zero ext}. 

We have $\mathrm{Ext}^1(\underline{S},X/\overline{X})=0$ if and only if at least one of the following two conditions is satisfied:
\begin{itemize}
    \item [(1)] $w_L=\alpha_L$ and $v_L=\beta_L$;
    \item [(2)] $w_R=\alpha_R$ and $v_R=\beta_R$.
\end{itemize}
In fact, (1) implies $b=d=\emptyset$ in Definition \ref{def_ext_gent_alg}, while (2) implies $a=c=\emptyset$.

In the rest of the proof we will repeatedly use Remark \ref{rmk_EC}. If $\mathrm{Ext}^1(\underline{S},X/\overline{X})\neq0$, then none of (1) and (2) is satisfied. Let $L$ be a successor closed subquiver of the diagram representing $X\oplus S$ which is not a successor closed subquiver of $Y$. We see from the description of the extension given in Theorem \ref{thm_A} that $L$ contains the diagram of $\underline{S}$ as a successor closed subquiver. Moreover, $L/\underline{S}$ is a successor closed subquiver of dimension $\textbf{dim}L-\textbf{dim}\underline{S}$ of $M=M(\alpha_L)\oplus M(\alpha_R) \oplus M(\beta_L)\oplus M(\beta_R)=\overline{X} \oplus S/\underline{S}$. See Figure \ref{non-zero ext}. This shows that every successor closed subquiver $L$ of the diagram of $X \oplus S$ that is not a subquiver of $Y$ gives rise to a successor closed subquiver of $M$ of dimension $\textbf{dim} L-\textbf{dim} \underline{S}$. Using Remark \ref{rmk_EC}, this completes the proof in this case. We observe that only two among $\alpha_L, \alpha_R, \beta_L, \beta_R$ can be non-zero. For example, $\alpha_L \neq 0$ and $ \alpha_R \neq 0$ imply that $\beta_L = \beta_R=0$.

\begin{figure}
    \centering
    \includegraphics[width=1\linewidth]{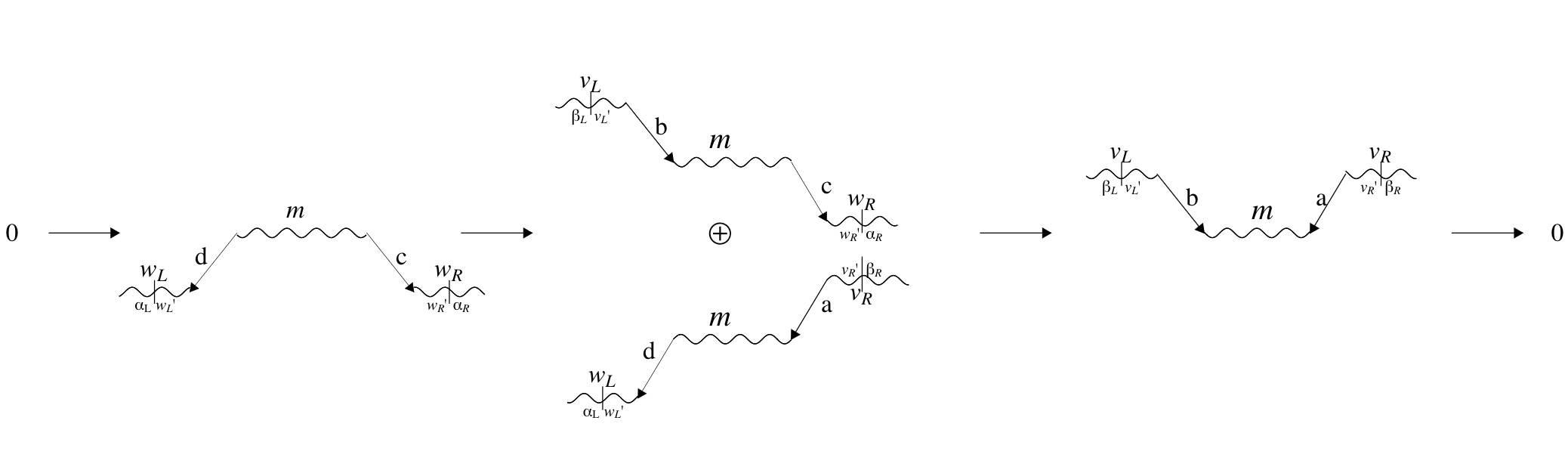}
    \caption{An overlap extension between $X=M(w)$ and $S=M(v)$ with $\mathrm{Ext}^1(\underline{S},X/\overline{X})\neq 0$.}
    \label{non-zero ext}
\end{figure}
\begin{figure}
    \centering
    \includegraphics[width=1\linewidth]{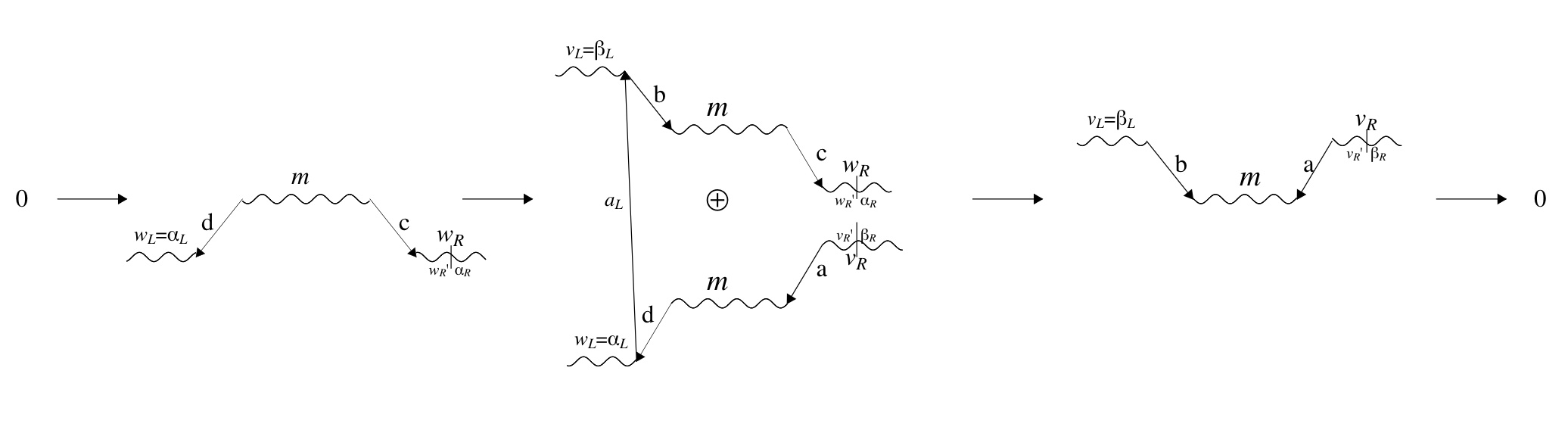}
    \caption{An overlap extension between $X=M(w)$ and $S=M(v)$ with $\mathrm{Ext}^1(\underline{S},X/\overline{X})= 0$ satisfying (1).}
    \label{zero ext 1}
\end{figure}
\begin{figure}
    \centering
    \includegraphics[width=1\linewidth]{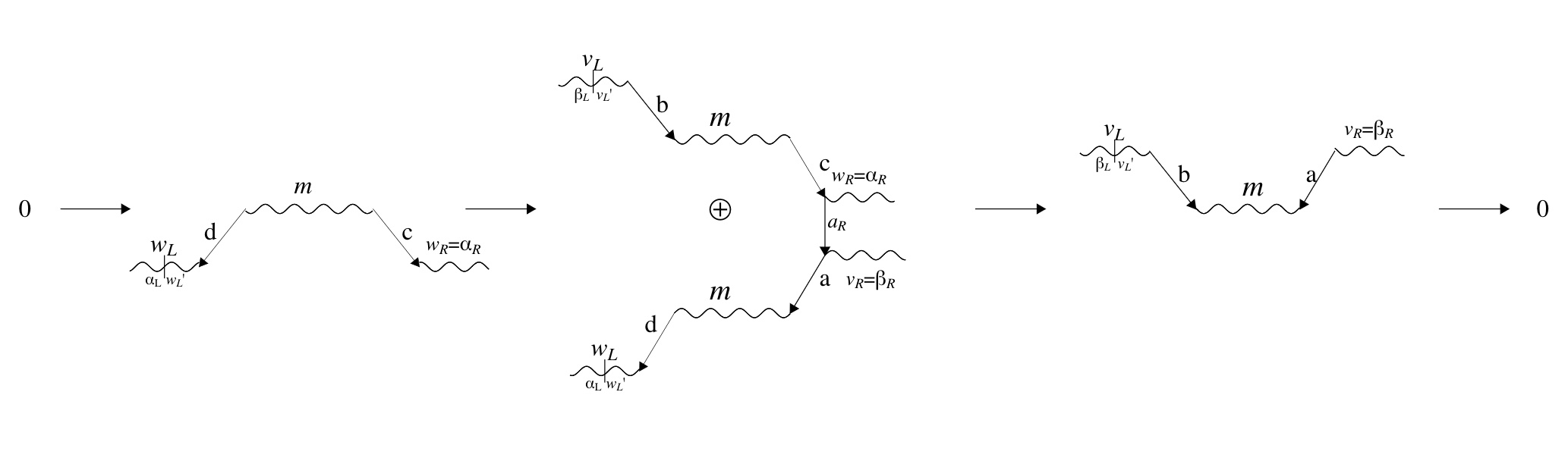}
    \caption{An overlap extension between $X=M(w)$ and $S=M(v)$ with $\mathrm{Ext}^1(\underline{S},X/\overline{X})= 0$ satisfying (2).}
    \label{zero ext 2}
\end{figure}
\begin{figure}
    \centering
    \includegraphics[width=1\linewidth]{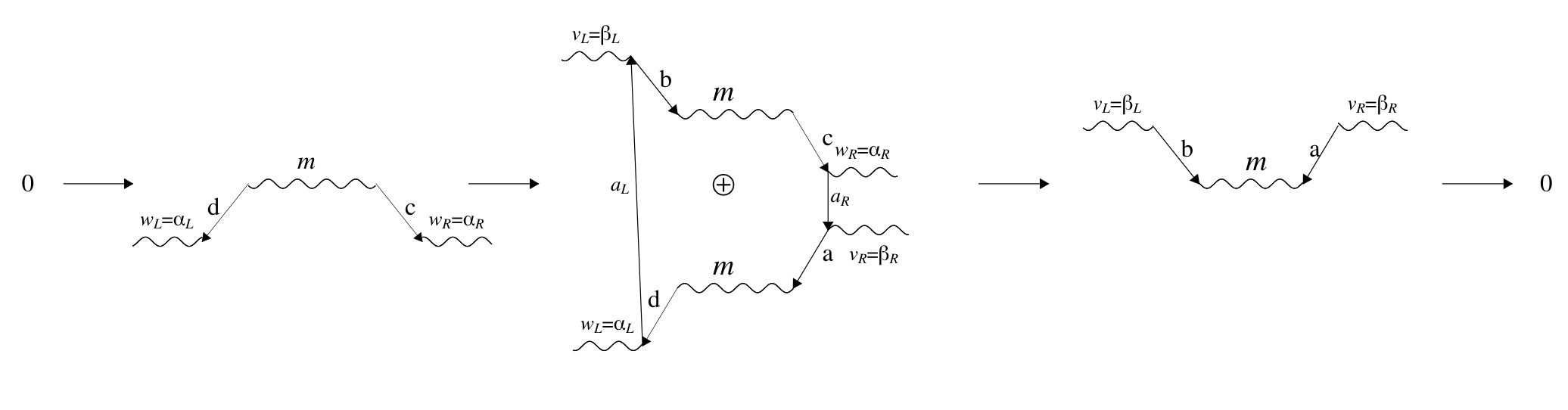}
    \caption{An overlap extension between $X=M(w)$ and $S=M(v)$ with $\mathrm{Ext}^1(\underline{S},X/\overline{X})=0$ satisfying (1) and (2).}
    \label{zero ext 3}
\end{figure}

Assume now that $\mathrm{Ext}^1(\underline{S},X/\overline{X})=0$, and only (1) is satisfied. Let $L$ be a successor closed subquiver of $X\oplus S$ which is not a successor closed subquiver of $Y$. Then $L$ contains the diagram of $\underline{S}$ as a successor closed subquiver. Moreover, $L/\underline{S}$ is a successor closed subquiver of dimension $\textbf{dim}L-\textbf{dim}\underline{S}$ of $M=M(w_La_Lv_L^{-1})\oplus M(\alpha_R) \oplus M(\beta_R)$, where $0 \to M(v_L) \to M(w_La_Lv_L^{-1}) \to M(w_L) \to 0$ is the non-trivial arrow extension between $M(v_L)$ and $M(w_L)$ given by the arrow $a_L:t(w_L) \to t(v_L) \in Q_1'$. In fact, if $L/\underline{S}$ had $M(w_L)$ as a successor closed subquiver, then $L$ would be a successor closed subquiver of $Y$. See Figure \ref{zero ext 1}. We have $M <_{\mathrm{Ext}} M(v_L)\oplus M(\alpha_R) \oplus M(v_L)\oplus M(\beta_R)=\overline{X} \oplus S/\underline{S}$. By Theorem \ref{thm_A} there are no non-trivial extensions between the indecomposable summands of $M$, therefore $M$ is the $\mathrm{Ext}$-minimum extension between $S/\underline{S}$ and $\overline{X}$. We observe that, in this case, only one between $\alpha_R$ and $\beta_R$ can be non-zero.

The case only (2) is analogous with $a_R:s(w_R)\to s(w_R) \in Q_1'$ and $M=M(\alpha_L)\oplus M(\beta_L)\oplus M(w_Ra_Rv_R^{-1})$. See Figure \ref{zero ext 2}.

If both (1) and (2) are satisfied, and $L$ is a successor closed subquiver of $X\oplus S$ that is not a successor closed subquiver of $Y$, then $L$ contains $\underline{S}$ as a successor closed subquiver, and $L/\underline{S}$ is a successor closed subquiver of dimension $\textbf{dim}L-\textbf{dim}\underline{S}$ of $M=M(w_La_Lv_L^{-1})\oplus M(w_Ra_Rv_R^{-1})$, where $a_L:t(w_L) \to t(v_L) \in Q_1'$ $a_R:s(w_R)\to s(w_R) \in Q_1'$. In fact, if $L/\underline{S}$ had $M(w_L)$ or $M(w_R)$ as a successor closed subquiver, then $L$ would be a successor closed subquiver of $Y$. See Figure \ref{zero ext 3}. We have $M  <_{\mathrm{Ext}} M(w_L)\oplus M(v_L) \oplus M(w_Ra_Rv_R^{-1}) <_{\mathrm{Ext}} M(v_L)\oplus M(w_R) \oplus M(v_L)\oplus M(v_R)=\overline{X} \oplus S/\underline{S}$, and $M  <_{\mathrm{Ext}} M(w_La_Lv_L^{-1}) \oplus M(w_R) \oplus M(v_R) <_{\mathrm{Ext}} M(v_L)\oplus M(w_R) \oplus M(v_L)\oplus M(v_R)=\overline{X} \oplus S/\underline{S}$. By Theorem \ref{thm_A} there are no non-trivial extensions between the indecomposable summands of $M$, therefore $M$ is the $\mathrm{Ext}$-minimum extension between $S/\underline{S}$ and $\overline{X}$.

By Remarks \ref{rmk_rel_jac} and \ref{rmk_rel_ext}, if $A=kQ/I$ is the gentle algebra of a triangulation of an unpuntured surface with marked points, then $a_L, a_R \in Q_1$. Therefore, $B=A$.

Finally, we consider the case where $\xi$ is an arrow extension given by $a:s(w)\to t(w) \in Q_1$. Since $[S,X]^1=1$, the morphism from $X=M(w)$ to $\tau S=M(g_r(g_l(v))$ (resp. from $\tau^{-1}X=M(f_l(f_r(w)))$ to $S=M(v)$) associated with the maximal overlap in $X$ between $g_l(v)$ and $w$ (resp. in $S$ between $v$ and $f_r(w)$) does not factor through an injective (resp. projective) module.


It follows that $\overline{X}=M(\alpha_L)$ and $S/\underline{S}=M(\beta_R)$, where $\alpha_L$ is the zero string or an image substring of $w$ such that $s(w)=s(\alpha_L)$, and $\beta_R$ is the zero string or a factor substring of $v$ such that $t(v)=t(\beta_R)$. Moreover, $a$ also gives an arrow extension between $X/X^S$ and $\underline{S}$, so $\mathrm{Ext}^1(\underline{S},X/\overline{X})\neq0$. Let $L$ be a successor closed subquiver of $X\oplus S$ which is not a successor closed subquiver of $Y$, then $L$ contains $\underline{S}$ as a successor closed subquiver. Moreover, $L/\underline{S}$ is a successor closed subquiver of dimension $\textbf{dim}L-\textbf{dim}\underline{S}$ of $M=M(\alpha_L) \oplus M(\beta_R)=\overline{X} \oplus S/\underline{S}$. This concludes the proof.

\end{proof}

We have the following immediate corollary:
\begin{corollary}\label{cor1}
  Let $A=kQ/I$ be a gentle algebra. Let $X$, $S$ be $A$-modules such that $[S,X]^1=1$. Let $\xi \in \mathrm{Ext}^1(S,X)$ be a generating extension with middle term $Y$. Then
    \begin{equation}\label{eq_cor1}
        CC(X)CC(S)= CC(Y)\bold{x}^{\bold{g}_X+\bold{g}_S-\bold{g}_Y} + \bold{y}^{\textbf{dim} \underline{S}}CC(M)\bold{x}^{B\textbf{dim}\underline{S}+\bold{g}_X+\bold{g}_S-\bold{g}_M}.
    \end{equation} 
\end{corollary}
\begin{remark}\label{rmk_spec}
    Specializing at $x_1=\cdots=x_n=1$, where $n=|Q_0|$, we get
     \begin{equation}\label{eq_rmk_cor1}
        F_XF_S= F_Y + \bold{y}^{\textbf{dim} \underline{S}}F_M.
    \end{equation} 
\end{remark}

\begin{examples}\label{ex_dec_thm}
    \begin{itemize}
        \item [(i)] Let $A$ be the algebra of Example \ref{gent_a_exs} (i). We consider the generating extension 
\begin{center}
    $\xi: 0 \to X=\begin{smallmatrix}
       2\\1
    \end{smallmatrix}=M(b)\to Y=\begin{smallmatrix}
        32\\1
    \end{smallmatrix}=M(bc^{-1})\to S=\begin{smallmatrix}
        3
    \end{smallmatrix}=M(1_3)\to 0$.
\end{center}

In this case, $\mathrm{Hom}(X,\tau S)=\mathrm{Hom}(\begin{smallmatrix}
       2\\1
    \end{smallmatrix},\begin{smallmatrix}
       32\\21
    \end{smallmatrix})=\mathrm{Hom}(1_2^-b,bc^{-1}a1_2^-)$ is of dimension 2, since there is the irreducible injective morphism $f$ that corresponds to the maximal overlap $b$ in $X$ between $bc^{-1}a1_2^-$ and $1_2^-b$, and there is also the morphism $g$ with kernel $\begin{smallmatrix}
    1
\end{smallmatrix}$ and image $\begin{smallmatrix}
    2
\end{smallmatrix}$ corresponding to the overlap $1_2^-$. As in the proof of Theorem \ref{dec_qg_gent_alg}, $f$ does not factor through an injective module, while $g$ factors through $I(1)=\begin{smallmatrix}
        32\\1
    \end{smallmatrix}$. Then $\overline{X}=\mathrm{ker}(f)=0$. 
    On the other hand, $\underline{S}=\mathrm{im}(\tau^{-1}X \to S)=\mathrm{im}(\begin{smallmatrix}
        3\\2
    \end{smallmatrix}\to \begin{smallmatrix}
        3 \end{smallmatrix})=\begin{smallmatrix}
        3
    \end{smallmatrix}$. So $M=0$.
    
    We have that $\mathrm{Ext}^1(\underline{S},X/\overline{X})=\mathrm{Ext}^1(\begin{smallmatrix}
        3
    \end{smallmatrix},\begin{smallmatrix}
        2\\1
    \end{smallmatrix})\neq 0$, and $F_{\overline{X}\oplus S/\underline{S}}=1$. In fact, one can compute that
    \begin{center}
        $F_{\begin{smallmatrix}
        2\\1
    \end{smallmatrix}}F_{\begin{smallmatrix}
        3
    \end{smallmatrix}}=F_{\begin{smallmatrix}
        32\\1
    \end{smallmatrix} }+y_3$.
    \end{center}
    \item [(ii)] Let $A$ be the algebra of Example \ref{ex_band}. We consider the generating extension
    \begin{center}
    $\xi: 0 \to X=M(a^{-1}ead^{-1}fd,1,\lambda)\to Y=M(cead^{-1}fda^{-1}ead^{-1})\to S=M(cead^{-1})\to 0$.
\end{center}
We have $\overline{X}=M(d)=\begin{smallmatrix}
    3\\1
\end{smallmatrix}$, $\underline{S}=M(ead^{-1})$, $S/\underline{S}=M(1_4)=\begin{smallmatrix}
    4
\end{smallmatrix}$. Moreover, with the notation of the proof of Theorem \ref{dec_qg_gent_alg}, $a_L=b \in Q_1$ gives the arrow extension
\begin{center}
    $0 \to S/\underline{S}=M(1_4)=\begin{smallmatrix}
        4
    \end{smallmatrix} \to M=M(db)=\begin{smallmatrix}
        3\\1\\4
    \end{smallmatrix} \to \overline{X}=M(b)=\begin{smallmatrix}
        3\\1
    \end{smallmatrix} \to 0$.
\end{center}
By Theorem \ref{dec_qg_gent_alg},
\begin{center}
    $F_XF_S=F_Y + y_1y_2^2y_3F_{\begin{smallmatrix}
        3\\1\\4
    \end{smallmatrix}}$.
\end{center}
See Example \ref{ex_F_pol} (ii) for the computation of $F_X$. In this case $B=A$, since $a_L=b \in Q_1$.
\item [(iii)] Let $A=kQ/I$ be the gentle algebra given by the quiver 
\begin{center}
    $\begin{tikzcd}
2 &                                    & 4 \arrow[ld, "c"']                 \\
  & 1 \arrow[lu, "a"'] \arrow[ld, "b"] &                                    \\
3 &                                    & 5 \arrow[lu, "d"] \arrow[uu, "e"']
\end{tikzcd}$
\end{center}
and $I=\langle ca,bd,ec\rangle$. We consider the following generating extension between $X=\begin{smallmatrix}
        5\\14\\2\hspace{0.1cm}
    \end{smallmatrix}=M(a^{-1}d^{-1}e)$ and $S=\begin{smallmatrix}
        45\\1
    \end{smallmatrix}=M(cd^{-1})$:
\begin{center}
    $\xi:0 \to X=\begin{smallmatrix}
        5\\14\\2\hspace{0.1cm}
    \end{smallmatrix}\to Y=\begin{smallmatrix}
        45\\14
    \end{smallmatrix}\oplus \begin{smallmatrix}
        5\\1\\2
    \end{smallmatrix}\to S=\begin{smallmatrix}
        45\\1
    \end{smallmatrix}\to 0$.
\end{center}
We have $\overline{X}=\begin{smallmatrix}
    2
\end{smallmatrix}$, $\underline{S}=\begin{smallmatrix}
    5\\1
\end{smallmatrix}=M(d)$, $S/\underline{S}=\begin{smallmatrix}
    4
\end{smallmatrix}$. In this case, $\mathrm{Ext}^1(\underline{S},X/\underline{X})=\mathrm{Ext}^1(\begin{smallmatrix}
    5\\1
\end{smallmatrix},\begin{smallmatrix}
    5\\1
\end{smallmatrix})=0$. Using the notation of the proof of Theorem \ref{dec_qg_gent_alg}, we have to consider an additional arrow $a_L:2 \to 4$. In this case, the algebra $B=kQ'/I$ that we obtain is still gentle. 
\begin{center}
    \hspace{-2cm}$Q': \begin{tikzcd}
      & 2 \arrow[rr, "a_L"] &                                    & 4 \arrow[ld, "c"']                 \\
&                     & 1 \arrow[lu, "a"'] \arrow[ld, "b"] &                                    \\
      & 3                   &                                    & 5 \arrow[lu, "d"] \arrow[uu, "e"']
\end{tikzcd}$
\end{center}
In $B$ we have the following arrow extension between $S/\underline{S}$ and $\overline{X}$:
\begin{center}
    $0 \to S/\underline{S}=\begin{smallmatrix}
        4
    \end{smallmatrix} \to M=\begin{smallmatrix}
        2\\4
    \end{smallmatrix}\to \overline{X}=\begin{smallmatrix}
        2
    \end{smallmatrix} \to 0$.
\end{center}
Therefore,
\begin{center}
    $F_{\begin{smallmatrix}
        5\\14\\2\hspace{0.1cm}
    \end{smallmatrix}}F_{\begin{smallmatrix}
        45\\1
    \end{smallmatrix}}=F_{\begin{smallmatrix}
        45\\14
    \end{smallmatrix}\oplus \begin{smallmatrix}
        5\\1\\2
    \end{smallmatrix}}+y_1y_5F_{\begin{smallmatrix}
        2\\4
    \end{smallmatrix}}$.
\end{center}
On the other hand, if we consider the generating extension between $X=\begin{smallmatrix}
    1\\3
\end{smallmatrix}=M(b)$ and $S=\begin{smallmatrix}
    45\\1
\end{smallmatrix}=M(d^{-1}e)$ given by
\begin{center}
    $\xi:0 \to X=\begin{smallmatrix}
    1\\3
\end{smallmatrix} \to Y= \begin{smallmatrix}
    4\\1\\3
\end{smallmatrix}\oplus \begin{smallmatrix}
    5\\1
\end{smallmatrix}\to S= \begin{smallmatrix}
    45\\1
\end{smallmatrix} \to 0$.
\end{center}
We have $\overline{X}=\begin{smallmatrix}
    3
\end{smallmatrix}$, $\underline{S}=\begin{smallmatrix}
    4\\1
\end{smallmatrix}=M(c)$, $S/\underline{S}=\begin{smallmatrix}
    5
\end{smallmatrix}$, and $\mathrm{Ext}^1(\underline{S},X/\underline{X})=\mathrm{Ext}^1(\begin{smallmatrix}
    4\\1
\end{smallmatrix},\begin{smallmatrix}
    1
\end{smallmatrix})=0$. We have to add an arrow $a_R:3 \to 5$. The algebra $B=kQ'/I$ that we obtain is no longer gentle; it is not even a string algebra.

\begin{center}
   \hspace{-2cm}$Q': \begin{tikzcd}
2                    &                                    & 4 \arrow[ld, "c"']                 \\
                     & 1 \arrow[lu, "a"'] \arrow[ld, "b"] &                                    \\
3 \arrow[rr, "a_R"'] &                                    & 5 \arrow[lu, "d"] \arrow[uu, "e"']
\end{tikzcd}$
\end{center}

In $B$, $M=\begin{smallmatrix}
        3\\5
    \end{smallmatrix} <_{\mathrm{Ext}}\overline{X}\oplus S/\underline{S}=\begin{smallmatrix}
        3
    \end{smallmatrix}\oplus \begin{smallmatrix}
        5
    \end{smallmatrix}$. Therefore,
\begin{center}
    $F_{\begin{smallmatrix}
        1\\3\hspace{0.1cm}
    \end{smallmatrix}}F_{\begin{smallmatrix}
        45\\1
    \end{smallmatrix}}=F_{\begin{smallmatrix}
        4\\1\\3
    \end{smallmatrix}\oplus \begin{smallmatrix}
        5\\1
    \end{smallmatrix}}+y_1y_4F_{\begin{smallmatrix}
        3\\5
    \end{smallmatrix}}$.
\end{center}
    \end{itemize}
\end{examples}

\begin{theorem}\label{thm_ca}
   Let $A=kQ/I$ be the gentle algebra of a triangulation $T$ of an unpunctured marked surface $(S,M)$. Let $X$, $S$ be rigid and indecomposable $A$-modules such that $[S,X]^1=1$. Let $\xi \in \mathrm{Ext}^1(S,X)$ be a generating extension with middle term $Y$. Then \ref{eq_cor1} is an exchange relation between the cluster variables $CC(X)$ and $CC(S)$ for the cluster algebra $\mathcal{A}_\bullet(Q)$ with principal coefficients in the initial seed whose exchange matrix is $B(Q)$.
\end{theorem}
\begin{proof}
First, we observe that $X$ and $S$ are string modules (cf. Example \ref{nonrigid_band}). It follows that $\underline{S}$ is also a string module, so $\underline{S}=M(z)$, and $B \textbf{dim}\underline{S}=(x_i)_{i \in Q_0}$ is the sum of the columns of $B$ indexed by the vertices that appear in $z$ counted with their multiplicities. For simplicity, assume that all the vertices of $z$ are distinct. Then
\begin{itemize}
    \item $x_i=-1$ if $i \to j \in Q_1$, and $j \not\in z$ or $j$ is the first vertex of $z$;
    \item $x_i=1$ if $i \gets j \in Q_1$, and $j \not\in z$ or $j$ is the last vertex of $z$;
    \item $x_i=-2$ if $j \gets i \to k \in z$;
    \item $x_i=2$ if $j \to i \gets k \in z$;
    \item $x_i=0$ if $j \to i \to k \in z$ or $j \gets i \gets k \in z$.
\end{itemize}
Therefore, by Theorem \ref{thm_g_vect} we have $B \textbf{dim}\underline{S}=-\bold{g}_{\underline{S}}-\bold{g}_{\tau \underline{S}}$. It follows that $B \textbf{dim}\underline{S}+\bold{g}_X + \bold{g}_S - \bold{g}_M \in \mathbb{Z}_{\geq 0}^{Q_0}$, and the only coordinates in $B \textbf{dim}\underline{S}+\bold{g}_X + \bold{g}_S - \bold{g}_M$ that might be non-zero correspond to vertices that are not in the strings of $X$ and $S$ nor in the strings of the indecomposable summands of $M$ and $Y$. The same happens for $\bold{g}_X+\bold{g}_S-\bold{g}_Y$. In fact, if $\xi$ is an arrow extension given by the arrow $a \in Q_1$, $\bold{g}_X+\bold{g}_S-\bold{g}_Y=\bold{e}_{s(\epsilon)}$, where $\epsilon \in Q_1$ is such that $v\epsilon^{-1}$ is a string but $wa^{-1}v\epsilon^{-1}$ is not a string. On the other hand, if $\xi$ is an overlap extension, $\bold{g}_X+\bold{g}_S-\bold{g}_Y=0$.

Let $\mathcal{C}_{(S,M)}$ be the cluster category of $(S,M)$ \cite{BZ,CSA}. Recall that an object $\Gamma \in \mathcal{C}_{(S,M)}$ is rigid if $\mathrm{Hom}_{\mathcal{C_{(S,M)}}}(\Gamma,\Gamma[1])=0$, where $[1]$ is the shift functor. In particular, rigid $A$-modules are rigid objects of $\mathcal{C}_{(S,M)}$. The cluster character is a bijection between rigid objects in $\mathcal{C}_{(S,M)}$ and cluster monomials in $\mathcal{A}_\bullet(Q)$. By assumption $X$ and $S$ are rigid. By Theorem \ref{thm_A} and Remark \ref{M_rigid}, we have that $X \oplus Y \oplus M$ and $S \oplus Y \oplus M$ are also rigid (as $A$-modules). Moreover, if $\bold{g}_X+ \bold{g}_S - \bold{g}_Y=(v_i)_{i \in Q_0}$ and $B \textbf{dim}\underline{S}+\bold{g}_X + \bold{g}_S - \bold{g}_M=(w_i)_{i \in Q_0}$, it follows from the previous observations that $\mathrm{Hom}(X,I(i))=\mathrm{Hom}(S,I(i))=\mathrm{Hom}(Y,I(i))=\mathrm{Hom}(M,I(i))=0$ for any $i \in Q_0$ such that $v_i > 0$ or $w_i > 0$. This concludes the proof.

\end{proof}

\begin{example}
    Let $A=kQ/I$ be the algebra of Example \ref{gent_a_exs} (ii). As observed in Example \ref{ex_alg_triang}, $A$ is the gentle algebra of the triangulation $T'$ of the octagon in Figure \ref{fig:ex_typeB}. We consider the following generating extension between the string modules $X=\begin{smallmatrix}
        3\\14
    \end{smallmatrix}=M(c^{-1}d)$ and $S=\begin{smallmatrix}
        25\\3
    \end{smallmatrix}=M(bf^{-1})$:

\begin{center}
    $\xi: 0 \to X=\begin{smallmatrix}
        3\\14
    \end{smallmatrix}\to Y=\begin{smallmatrix}
        5\\3\\1
    \end{smallmatrix}\oplus \begin{smallmatrix}
        2\\3\\4
    \end{smallmatrix} \to S=\begin{smallmatrix}
        25\\3
    \end{smallmatrix}\to 0$
\end{center}
 Then $\overline{X}=\begin{smallmatrix}
        1
    \end{smallmatrix}\oplus\begin{smallmatrix}
        4
    \end{smallmatrix}$, $\underline{S}=\begin{smallmatrix}
        3
    \end{smallmatrix}$, $S/\underline{S}=\begin{smallmatrix}
        2
    \end{smallmatrix}\oplus \begin{smallmatrix}
        5
    \end{smallmatrix}$, and $\mathrm{Ext}^1(\underline{S},X/\overline{X})=\mathrm{Ext}^1(\begin{smallmatrix}
        3
    \end{smallmatrix},\begin{smallmatrix}
        3
    \end{smallmatrix})=0$. The $\leq_{\mathrm{Ext}}$-minimum extension between $S/\underline{S}$ and $\overline{X}$ is $M=\begin{smallmatrix}
        1\\2
    \end{smallmatrix}\oplus \begin{smallmatrix}
        4\\5
    \end{smallmatrix}$. We observe that, with the notation of the proof of Theorem \ref{dec_qg_gent_alg}, $a_L=a:1 \to 2 \in Q_1$ and $a_R=e : 4 \to 5 \in Q_1$, so $B=A$. Moreover, $\bold{g}_X + \bold{g}_S - \bold{g}_Y=B\text{dim}\underline{S}+\bold{g}_X + \bold{g}_S-\bold{g}_M=\bold{0}$.

    Therefore, by Theorem \ref{thm_ca}, we have the following exchange relation in the cluster algebra $\mathcal{A}_\bullet(Q)$ of type $A_5$:
    \begin{center}
        $CC({\begin{smallmatrix}
        3\\14
    \end{smallmatrix}})CC({\begin{smallmatrix}
        25\\3
    \end{smallmatrix}})=CC({\begin{smallmatrix}
        5\\3\\1
    \end{smallmatrix}\oplus \begin{smallmatrix}
        2\\3\\4
    \end{smallmatrix} })+y_3CC({\begin{smallmatrix}
        1\\2
    \end{smallmatrix}\oplus\begin{smallmatrix}
        4\\5
    \end{smallmatrix}})$.
    \end{center}
   
\end{example}

\begin{remark}
Let $\Lambda$ be a string algebra which is a quotient of a gentle algebra $A$. Let $X,S$ be $\Lambda$-modules such that $\text{dim}\mathrm{Ext}_A^1(S,X)^1=1$, and let $\xi:0 \to X \to Y \to S \to 0$ be a generating extension such that $Y$ is also a module over $\Lambda$. Then the proof of Theorem \ref{dec_qg_gent_alg} applies also to this slightly more general setting.
\end{remark}
\begin{example}
Let $Q$ be the quiver 
\begin{center}
    $\begin{tikzcd}
8 \arrow[rd, "b_1"'] &                                         & 1 \arrow[ll, "a_1"'] \arrow[rd, "c_2"'] &                                         & 2 \arrow[ll, "b_2"'] \\
                     & 7 \arrow[ru, "c_1"'] \arrow[ld, "a_4"'] &                                         & 3 \arrow[ru, "a_2"'] \arrow[ld, "c_3"'] &                      \\
6 \arrow[rr, "b_4"'] &                                         & 5 \arrow[lu, "c_4"'] \arrow[rr, "a_3"'] &                                         & 4 \arrow[lu, "b_3"']
\end{tikzcd}$
\end{center}
Let $I=\langle a_ib_i,b_ic_i,c_ia_i\rangle_{i=1}^4$ and $J=\langle c_1c_2c_3,c_2c_3c_4,c_3c_4c_1,c_4c_1c_2\rangle$. Then $A=kQ/I$ is gentle, while $\Lambda=A/J$ is a non-gentle string algebra. In particular, $\Lambda$ is the Jacobian algebra of the quiver with potential associated with a triangulation of the octagon with one puncture (cf. \cite{LF}). We consider the following generating extension between $X=\begin{smallmatrix}
        \hspace{0.1cm}58\\47
    \end{smallmatrix}=M(a_3^{-1}c_4b_1^{-1})$ and $S=\begin{smallmatrix}
        \hspace{0.1cm}36\\25
    \end{smallmatrix}=M(a_2^{-1}c_3b_4^{-1})$:
\begin{center}
    $0 \to X=\begin{smallmatrix}
        \hspace{0.1cm}58\\47
    \end{smallmatrix} \to Y=\begin{smallmatrix}
        \hspace{-0.1cm}3\\258\\ \hspace{0.1cm}7
    \end{smallmatrix}\oplus \begin{smallmatrix}
        6\\5\\4
    \end{smallmatrix}\to S=\begin{smallmatrix}
        \hspace{0.1cm}36\\25
    \end{smallmatrix}\to 0$.
\end{center}

We have $\overline{X}=\begin{smallmatrix}
    4
\end{smallmatrix}\oplus \begin{smallmatrix}
    8\\7
\end{smallmatrix}$, $\underline{S}=\begin{smallmatrix}
    5
\end{smallmatrix}$, $S/\underline{S}=\begin{smallmatrix}
    2\\3
\end{smallmatrix}\oplus \begin{smallmatrix}
    6
\end{smallmatrix}$, $\mathrm{Ext}^1(\underline{S},X/\underline{X})=\mathrm{Ext}^1(\begin{smallmatrix}
    5
\end{smallmatrix},\begin{smallmatrix}
    5
\end{smallmatrix})=0$, and $M=\begin{smallmatrix}
    4\\3\\2
\end{smallmatrix} \oplus \begin{smallmatrix}
    8\\7\\6
\end{smallmatrix}$. In this case, using the notation of the proof of Theorem \ref{dec_qg_gent_alg}, $a_L=b_3:4 \to 3 \in Q_1$ and $a_R=a_4 : 7 \to 6 \in Q_1$, so $B=A$. Moreover, $\bold{g}_X + \bold{g}_S - \bold{g}_Y=B\text{dim}\underline{S}+\bold{g}_X + \bold{g}_S-\bold{g}_M=\bold{0}$.

Therefore, we have:
    \begin{center}
        $CC({\begin{smallmatrix}
         \hspace{0.1cm}58\\47
    \end{smallmatrix}})CC({\begin{smallmatrix}
         \hspace{0.1cm}36\\25
    \end{smallmatrix}})=CC(\begin{smallmatrix}
        \hspace{-0.1cm}3\\258\\ \hspace{0.1cm}7
    \end{smallmatrix}\oplus \begin{smallmatrix}
        6\\5\\4
    \end{smallmatrix})+y_5CC(\begin{smallmatrix}
    4\\3\\2
\end{smallmatrix} \oplus \begin{smallmatrix}
    8\\7\\6
\end{smallmatrix})$.
    \end{center}

We observe that this is an exchange relation in the cluster algebra of type $D_8$ with principal coefficients in the initial seed whose exchange matrix is $B(Q)$.
\end{example}
\section{Application to cluster algebras of type B}\label{section_B}
Let $n$ be a positive integer, and let $\mathbf{P}_{n+3}$ be the regular polygon with $n+3$ vertices. Clusters of cluster algebras of type $A_n$ are in bijection with triangulations of $\mathbf{P}_{n+3}$, while cluster variables $x_\gamma$ correspond to the diagonals $\gamma$ of the polygon \cite{CAII}, with the convention that $x_\gamma = 1$ if $\gamma$ is a boundary segment. 

Let $\Bar{T}=\{\tau_1,\dots,\tau_n\}$ be a triangulation of $\mathbf{P}_{n+3}$. The cluster algebra $\mathcal{A}_\bullet(\Bar{T})$ of type $A_n$ with principal coefficients in $\Bar{T}$ is defined as the cluster algebra with principal coefficients in the initial seed whose exchange matrix is the skew-symmetric matrix $(b_{ij})$ with entries in $\{-1,0,1\}$, such that $b_{ij}=1$ if and only if $\tau_i$ and $\tau_j$ are two sides of a triangle in $\Bar{T}$, and $\tau_i$ follows $\tau_j$ in counterclockwise order. 

Let $(a,b)$ be the diagonal of $\mathbf{P}_{n+3}$ which connects $a$ and $b$. The \emph{elementary lamination} associated with $(a,b)$ is the segment which begins at a point $a'$ of the polygon located near $a$ in the clockwise direction, and ends at a point $b'$ near $b$ in
the clockwise direction, as in Figure \ref{lamination}. Elementary laminations provide information about the coefficients in exchange relations in $\mathcal{A}_\bullet(\Bar{T})$, as showed in the following result:

 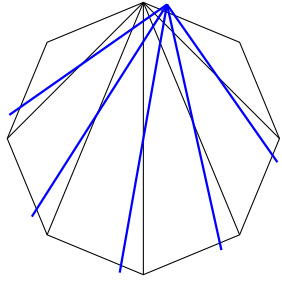
\begin{figure}
        \centering
\begin{tikzpicture}[scale=0.6]
     \draw (90:3cm) -- (135:3cm) -- (180:3cm) -- (225:3cm) -- (270:3cm) -- (315:3cm) -- (360:3cm) -- (45:3cm) -- cycle;
                    \draw (90:3cm) -- node[near end, below, xshift=-1mm, yshift=-1mm] {} (180:3cm);
                    \draw (90:3cm) -- node[near end,below,xshift=1.2mm] {} (225:3cm);
                    \draw (90:3cm) -- node[near end, below,xshift=1mm] {} (270:3cm);
                    \draw (90:3cm) -- node[near end, above, xshift=-1mm] {} (315:3cm);
                    \draw (90:3cm) -- node[near end, above, xshift=-1mm] {} (360:3cm);

                    \draw[blue, line width=0.3mm] (80:3cm) -- node[near end, above, xshift=-1mm, yshift=-1mm] {} (170:3cm);
                    \draw[blue, line width=0.3mm] (80:3cm) -- node[near end, above left,xshift=1.2mm] {} (215:3cm);
                    \draw[blue, line width=0.3mm] (80:3cm) -- node[near end, above left,xshift=1mm] {} (260:3cm);
                    \draw[blue, line width=0.3mm] (80:3cm) -- node[near end, below, xshift=-1.5mm] {} (305:3cm);
                    \draw[blue, line width=0.3mm] (80:3cm) -- node[near end, below] {} (350:3cm);

\end{tikzpicture}
        \caption{A triangulated octagon with the elementary lamination associated with each diagonal of the triangulation (in blue).}
        \label{lamination}
    \end{figure}

\begin{proposition}[{\cite[Proposition 17.3]{FT}}] \label{up:skein1}
Let $(a,b)$ and $(c,d)$ be two diagonals
which intersect each other. Then
\begin{equation} \label{u:skein-eq1}
x_{(a,b)} x_{(c,d)} = \bold{y}^{\bold{d}_{ac,bd}} x_{(a,d)} ~x_{(b,c)}
 + \bold{y}^{\bold{d}_{ad,bc}}x_{(a,c)} ~x_{(b,d)},
\end{equation}
where
$\bold{d}_{ac,bd}$ (resp., $\bold{d}_{ad,bc}$) is the vector whose $i$-th coordinate is 1 if the elementary lamination associated with $\tau_i$ crosses both $(a,c)$ and $(b,d)$ (resp., $(a,d)$ and $(b,c)$); 0 otherwise.
\end{proposition}

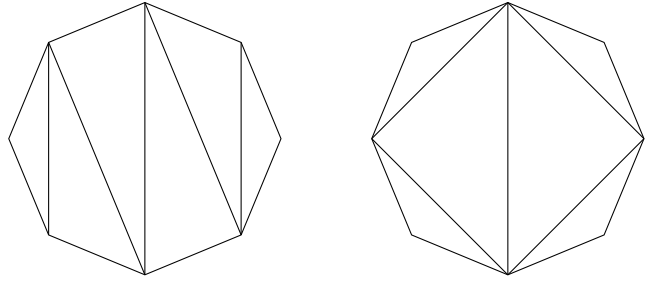
\begin{figure}[H]
                               \centering
                \begin{tikzpicture}[scale=0.6]
                    \draw (90:3cm) -- (135:3cm) -- (180:3cm) -- (225:3cm) -- (270:3cm) -- (315:3cm) -- (360:3cm) -- (45:3cm) -- cycle;
                    \draw (135:3cm) -- (225:3cm);;
                    \draw (45:3cm) -- (315:3cm);
                   \draw[] (270:3cm) -- node[midway, above left,xshift=1mm] {} (90:3cm);
                    \draw (135:3cm) -- node[midway, above left,xshift=1mm] {} (270:3cm);
                    \draw (315:3cm) -- node[midway, above left,xshift=1mm] {} (90:3cm);

                    \begin{scope}[xshift=8cm]
                      \draw (90:3cm) -- (135:3cm) -- (180:3cm) -- (225:3cm) -- (270:3cm) -- (315:3cm) -- (360:3cm) -- (45:3cm) -- cycle;
                    \draw (360:3cm) -- (270:3cm);
                    \draw (180:3cm) -- node[midway, above left,xshift=1mm] {} (90:3cm);
                    \draw (90:3cm) -- node[midway, above left,xshift=1mm] {} (270:3cm);
                     \draw (360:3cm) -- node[midway, above left,xshift=1mm] {} (90:3cm);
                    \draw (180:3cm) -- node[midway, above left,xshift=1mm] {} (270:3cm);

                    \end{scope}
                \end{tikzpicture}
                               \caption{Two $\theta$-invariant triangulations of $\mathbf{P}_8$.}
                               \label{theta-inv-triang}
                           \end{figure}

On the other hand, let $\mathbf{P}_{2n+2}$ be the regular polygon with $2n+2$ vertices, and let $\theta$ be the rotation of $180^\circ$. Clusters of cluster algebras of type $B_n$  are in bijection with $\theta$-invariant triangulations of $\mathbf{P}_{2n+2}$, see Figure \ref{theta-inv-triang} for an example, and cluster variables $x_{ab}$ correspond to the orbits $[a,b]$ of the action of $\theta$ on the diagonals $(a,b)$ of the polygon \cite{CAII}, with the convention that $x_{ab}=1$ if $[a,b]$ lies on the boundary of $\mathbf{P}_{2n+2}$. 

Each $\theta$-invariant triangulation $T$ of $\mathbf{P}_{2n+2}$ has exactly one diameter $d$. After choosing an orientation of $d$, one can define the cluster algebra $\mathcal{A}_\bullet^B(T)$ of type $B_n$ with principal coefficients in $T$ (\cite{ciliberti2024}, Definition 1.3). It turns out that cluster variables in $\mathcal{A}_\bullet ^B(T)$ are related to the cluster variables of type $A_n$ thanks to the following simple operation on the diagonals of the polygon:

\begin{definition}[{\cite[Definition 3.1]{ciliberti2024}}]\label{def_restriction}
    Let $\mathcal{D}$ be a set of diagonals of $\mathbf{P}_{2n+2}$. The \emph{restriction of $\mathcal{D}$}, denoted by $\text{Res}(\mathcal{D})$, is the set of diagonals of $\mathbf{P}_{n+3}$ obtained from those of $\mathcal{D}$ identifying all the vertices which lie on the right of $d$.
\end{definition}

Let $T=\{ \tau_1, \dots, \tau_n=d, \dots, \tau_{2n-1} \}$ be a $\theta$-invariant triangulation of $\mathbf{P}_{2n+2}$ with oriented diameter $d$. Let $x_{ab}$ be the cluster variable of type $B_n$ in $\mathcal{A}_\bullet^B(T)$ corresponding to the $\theta$-orbit $[a,b]$ of the diagonal $(a,b)$. The $F$-polynomial and the $\bold{g}$-vector of $x_{ab}$ are denoted by $F_{ab}$ and $\bold{g}_{ab}$. Moreover, for a diagonal $\gamma$, $F_\gamma$ and $\bold{g}_\gamma$ denote the $F$-polynomial and the $\bold{g}$-vector of $x_\gamma \in \mathcal{A}_\bullet(\Bar{T})$, where $\Bar{T}=\text{Res}(T)=\{\tau_1,\dots,\tau_n\}$ is the triangulation of $\mathbf{P}_{n+3}$ obtained after applying the restriction to $T$.

\begin{theorem}[{\cite[Theorem 3.7]{ciliberti2024}}]\label{t1} Let $[a,b] \not \subset T$ be a $\theta$-orbit. Let $D=\text{diag}(1,\dots,1,2)\in \mathbb{Z}^{n\times n}$.
\begin{itemize}
                           \item[(i)] If $\text{Res}([a,b])$ contains only one diagonal $\gamma$ (as in Figure \ref{type B 1}), then 
                           \begin{align}\label{all1}
                               F_{ab}&=F_\gamma,\\
                               \bold{g}_{ab}&=\begin{cases}
          \text{$D\bold{g}_\gamma$ if $\gamma$ does not cross $\tau_n=d$;}\\
          \text{$D\bold{g}_\gamma+\bold{e}_n$ if $\gamma$ crosses $\tau_n=d$}.
      \end{cases}
                           \end{align}
                           \item[(ii)] Otherwise, $\text{Res}([a,b])=\{ \gamma_1, \gamma_2\}$ (as in Figure \ref{type B}), and
                           \begin{align}\label{all2}
                               F_{ab}&=F_{\gamma_1}F_{\gamma_2}- \bold{y}^{\bold{d}_{\gamma_1,\gamma_2}}F_{(a,\theta(b))},\\
                               \bold{g}_{ab}&=D(\bold{g}_{\gamma_1}+\bold{g}_{\gamma_2}+\bold{e}_n).
                           \end{align}

 where $\bold{d}_{\gamma_1,\gamma_2} \in \{0,1\}^n$ is such that $(\bold{d}_{\gamma_1,\gamma_2})_i=1$ if and only if the elementary lamination associated with $\tau_i$ crosses both $\gamma_1$ and $\gamma_2$, $i=1,\dots,n$, and $\bold{e}_n$ is the $n$-th vector of the canonical basis of $\mathbb{Z}^n$.
                             \end{itemize}
   \begin{figure}
                               \centering
                \begin{tikzpicture}[scale=0.6]
                    \draw (90:3cm) -- (120:3cm) -- (150:3cm) -- (180:3cm) -- (210:3cm) -- (240:3cm) -- (270:3cm) -- (300:3cm) -- (330:3cm) -- (360:3cm) -- (30:3cm) -- (60:3cm) --  cycle;
                  
                    \draw[-{Latex[length=2mm]}] (270:3cm) -- node[midway, above left,xshift=1mm] {} (90:3cm);
                  
                     \draw[red, line width=0.3mm] (150:3cm) -- (330:3cm);

                    \node at (330:3cm) [right] {$\Bar{a}$};

\node at (150:3cm) [left] {$a$};

                    \begin{scope}[xshift=8cm]
                      \draw (90:3cm) -- (120:3cm) -- (150:3cm) -- (180:3cm) -- (210:3cm) -- (240:3cm) -- (270:3cm) -- (360:3cm)  --  cycle;
                    \draw (90:3cm) -- node[midway, above left,xshift=1mm] {} (270:3cm);
                    \draw[red, line width=0.3mm] (150:3cm) -- node[midway, above right, xshift=-2mm, yshift=-1mm] {$\textcolor{black}{\gamma}$} (360:3cm); 
                  
                    \node at (360:3cm) [right] {$\ast$};

\node at (150:3cm) [left] {$a$};
                      
                    \end{scope}
                    \begin{scope}[yshift=-9cm]
                         \draw (90:3cm) -- (120:3cm) -- (150:3cm) -- (180:3cm) -- (210:3cm) -- (240:3cm) -- (270:3cm) -- (300:3cm) -- (330:3cm) -- (360:3cm) -- (30:3cm) -- (60:3cm) --  cycle;
                  
                    \draw[-{Latex[length=2mm]}] (270:3cm) -- node[midway, above left,xshift=1mm] {} (90:3cm);
                  
                     \draw[red, line width=0.3mm] (150:3cm) -- (270:3cm); 
                   
 \draw[red, line width=0.3mm] (90:3cm) -- (330:3cm);

                    \node at (90:3cm) [above] {$\Bar{b}$};
                    \node at (330:3cm) [right] {$\Bar{a}$};

\node at (150:3cm) [left] {$a$};
                    \node at (270:3cm) [below] {$b$};

                    \begin{scope}[xshift=8cm]
                      \draw (90:3cm) -- (120:3cm) -- (150:3cm) -- (180:3cm) -- (210:3cm) -- (240:3cm) -- (270:3cm) -- (360:3cm)  --  cycle;
                    \draw (90:3cm) -- node[midway, above left,xshift=1mm] {} (270:3cm);
                    \draw[red, line width=0.3mm] (150:3cm) -- node[midway, above right, xshift=-2mm, yshift=-1mm] {$\textcolor{black}{\gamma}$} (270:3cm); 
                  
                    \node at (360:3cm) [right] {$\ast$};
 \node at (270:3cm) [below] {$b$};
\node at (150:3cm) [left] {$a$};
                      
                    \end{scope}
                    \end{scope}
                \end{tikzpicture}
                               \caption{On the left, two $\theta$-orbits $[a,\Bar{a}]$ and $[a,b]$. On the right, their restrictions.}
                               \label{type B 1}
                           \end{figure}
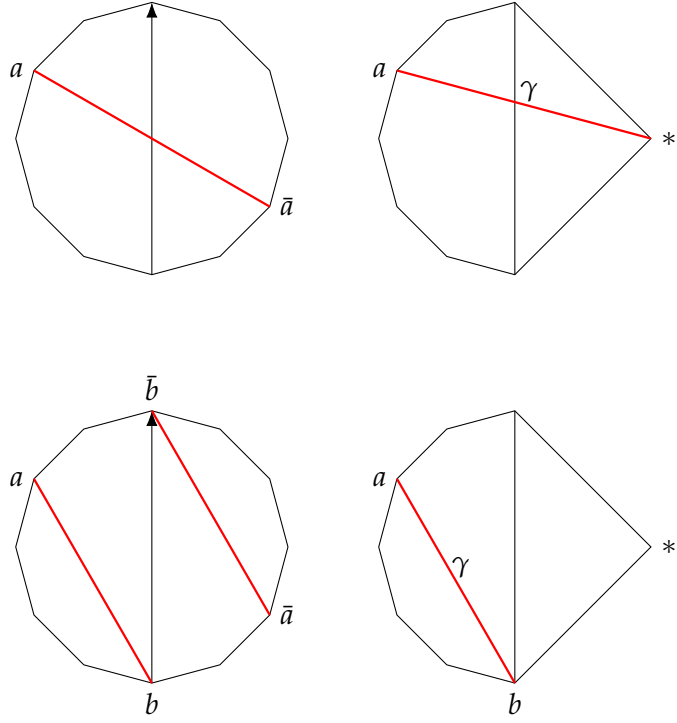  

  \begin{figure}
                               \centering
                \begin{tikzpicture}[scale=0.6]
                    \draw (90:3cm) -- (120:3cm) -- (150:3cm) -- (180:3cm) -- (210:3cm) -- (240:3cm) -- (270:3cm) -- (300:3cm) -- (330:3cm) -- (360:3cm) -- (30:3cm) -- (60:3cm) --  cycle;
                    \draw[-{Latex[length=2mm]}] (270:3cm) -- node[midway, above left,xshift=1mm] {} (90:3cm);
                     
                     \draw[red, line width=0.3mm] (150:3cm) -- (30:3cm); 
                    \draw[red, line width=0.3mm] (210:3cm) -- (330:3cm);

                    \node at (210:3cm) [left] {$\Bar{b}$};
                    \node at (330:3cm) [right] {$\Bar{a}$};

\node at (150:3cm) [left] {$a$};
                    \node at (30:3cm) [right] {$b$};

                    \begin{scope}[xshift=8cm]
                      \draw (90:3cm) -- (120:3cm) -- (150:3cm) -- (180:3cm) -- (210:3cm) -- (240:3cm) -- (270:3cm) -- (360:3cm)  --  cycle;
                    \draw (90:3cm) -- node[midway, above left,xshift=1mm] {} (270:3cm);
                    \draw[red, line width=0.3mm] (150:3cm) -- node[midway, above right, xshift=-2mm, yshift=-1mm] {$\textcolor{black}{\gamma_1}$} (360:3cm); 
                     \draw[blue, line width=0.3mm] (150:3cm) -- (210:3cm);
\draw[red, line width=0.3mm] (210:3cm) -- node[midway, above right, xshift=-1mm] {$\textcolor{black}{\gamma_2}$} (360:3cm);
                  
                    \node at (210:3cm) [left] {$\Bar{b}$};
                    \node at (360:3cm) [right] {$\ast$};

\node at (150:3cm) [left] {$a$};
                      
                    \end{scope}
                \end{tikzpicture}
                               \caption{On the left, a $\theta$-orbit $[a,b]$. On the right, its restriction in red and the diagonal $(a,\Bar{b})$ in blue.}
                               \label{type B}
                           \end{figure}
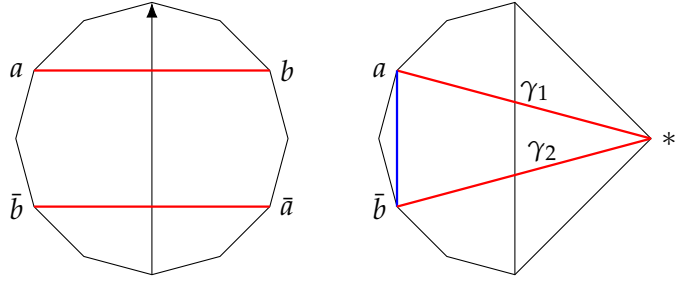

\end{theorem}
Furthermore, in \cite{ciliberti2024}, each cluster variable of type $B_n$ in $\mathcal{A}_\bullet^B(T)$ is associated with an indecomposable symmetric module over a symmetric algebra constructed from $T$. We briefly recall this construction; for more details refer to \cite[Section 4]{ciliberti2024}. 

First, we introduce some basic notions of symmetric representation theory developed by Derksen and Weyman in \cite{DW}, and by Boos and Cerulli Irelli in \cite{boos2021degenerations}, in order to set up the notation.

\begin{definition}\label{def_symm_quiver}
    A \emph{symmetric quiver} is a pair $(Q,\sigma)$, where $Q$ is a finite quiver and $\sigma$ is an involution of $Q_0$ and of $Q_1$ which reverses the orientation of arrows. 
\end{definition}

\begin{example}
    Let $Q = 1 \xrightarrow[]{a} 2 \xrightarrow[]{b} 3$ and $Q'= 1 \xrightarrow[]{a} 2 \xleftarrow[]{b} 3$ be two quivers of type $A_3$. Then $Q$ is symmetric, with the involution $\sigma$ given by $\sigma(1)=3$, $\sigma(2)=2$ and $\sigma(a)=b$, while $Q'$ is not symmetric, i.e., it cannot be endowed with the structure of a symmetric quiver.
\end{example}

\begin{definition}
    Let $(Q,\sigma)$ be a symmetric quiver. Let $I \subset kQ$ be an admissible ideal such that $\sigma(I)=I$. Then $A=kQ/I$ is called a \emph{symmetric quiver algebra}.
\end{definition}

\begin{definition}
    A \emph{symmetric module} over a symmetric algebra $A=kQ/I$ is a triple 
    $(V_i,\phi_a, \langle \cdot, \cdot \rangle)$, where $(V_i,\phi_a)$ is a $A$-module, $\langle \cdot, \cdot \rangle$ is a non-degenerate symmetric or skew-symmetric scalar product on $V=\displaystyle\bigoplus_{i \in Q_0}V_i$ such that its restriction to $V_i \times V_j$ is 0 if $j \neq \sigma(i)$, and $\langle \phi_a(v), w \rangle + \langle v, \phi_{\sigma(a)}(w) \rangle=0$, for every $a : i \to j \in Q_1$, $v \in V_i$, $w \in V_{\sigma(j)}$. If $\langle \cdot, \cdot \rangle$ is symmetric (resp. skew-symmetric), $(V_i,\phi_a, \langle \cdot, \cdot \rangle)$ is called $orthogonal$ (respectively, $symplectic$).
\end{definition}

\begin{definition}
    Let $L=(V_i,\phi_a)$ be a module over a symmetric algebra $A=kQ/I$. The $twisted$ $dual$ of $L$ is the $A$-module $\nabla L = (\nabla V_i, \nabla \phi_a)$, where $\nabla V_i=V_{\sigma(i)}^\ast$ and $\nabla \phi_a = - \phi_{\sigma (a)}^\ast$ ($\ast$ denotes the linear dual).
\end{definition}
\begin{remark}
 The twisted dual is a contravariant exact endofunctor on $\mathrm{mod}A$. Moreover, if $L$ is symmetric, the scalar product $\langle \cdot, \cdot \rangle$ induces an isomorphism from $V=\displaystyle\bigoplus_{i \in Q_0}V_i$ to $\nabla V=\displaystyle\bigoplus_{i \in Q_0}\nabla V_i$. 
\end{remark}


The following result shows that every indecomposable symmetric module is uniquely determined by the $\nabla$-orbit of an ordinary indecomposable module:

\begin{lemma}[{\cite[Lemma 2.10]{boos2021degenerations}}]\label{ind_symm}
 Let $N$ be an indecomposable symmetric module of a symmetric algebra $A$. Then, one and only one of the following three cases can occur: 
 \begin{itemize}
     \item [(I)] $N$ is indecomposable as a $A$-module; in this case, $N$ is called of type (I), for “indecomposable”;
     \item [(S)] there exists an indecomposable $A$-module $L$ such that $N=L\oplus \nabla L$ and $L\ncong \nabla L$; in this case, $N$ is called of type (S), for “split”;
     \item [(R)] there exists an indecomposable $A$-module $L$ such that $N=L\oplus \nabla L$ and $L\cong \nabla L$; in this case, $N$ is called of type (R) for “ramified”.
 \end{itemize}
\end{lemma}

\begin{lemma}\label{lemma_M}
If $A$ is a symmetric algebra and $L$ is an $A$-module such that $[\nabla L, L]^1=1$, then $M$ in \ref{dec_qg_gent_alg} is a $\nabla$-invariant $A$-moodule. 
\end{lemma}
\begin{proof}
    By definition, $\overline{L}= \text{ker}(L \to \tau \nabla L)$, and $\underline{\nabla L}= \text{im}(\tau^{-1} L \to \nabla L)$. Since $\nabla\tau = \tau^{-1} \nabla$ (cf. \cite{DW}, Proposition 3.4), $\nabla (\overline{L})=\text{coker}(\tau^{-1} L \to \nabla L)=\nabla L / \underline{\nabla L}$, $\overline{L}\oplus \nabla L / \underline{\nabla L}$ is $\nabla$-invariant. Consider the short exact sequence 
    \begin{center}
        $0 \to \nabla L / \underline{\nabla L} \to M \to \overline{L}\to 0$.
    \end{center}
    Since $\nabla$ is a contravariant exact endofunctor, the sequence
    \begin{center}
        $0 \to \nabla L / \underline{\nabla L} \to \nabla M \to \overline{L}\to 0$
    \end{center}
    is also exact. Assume $\nabla M \neq M$; then $\nabla M >_{\mathrm{Ext}} M$, by the minimality of $M$. So $\nabla M= M_1 \oplus M_2$ and there exists a non-split short exact sequence
    \begin{center}
        $0 \to M_1 \to M \to M_2 \to 0$.
    \end{center}
    Applying $\nabla$ to it, we have that also
    \begin{center}
        $0 \to \nabla M_2 \to \nabla M \to \nabla M_1 \to 0$
    \end{center}
    is exact and non-split. Since $M=\nabla M_1 \oplus \nabla M_2$, this contradicts the minimality of $M$.
\end{proof} 

Let $T=\{\tau_1,\dots,\tau_{2n-1}\}$ be a triangulation of $\mathbf{P}_{2n+2}$. Naturally associated with $T$ is a quiver $Q(T)$, where there is an arrow from vertex $j$ to vertex $i$ if and only if $\tau_i$ and $\tau_j$ are sides of a triangle in $\Bar{T}$, and $\tau_i$ is counterclockwise from $\tau_j$, and an ideal $I(T)$ generated by all paths $i \to j \to k$ such that there exists an arrow $k \to i$. See \cite{FST,LF} for more details. Moreover, with each diagonal $\alpha \not\subset T$ is associated an $A=kQ(T)/I(T)$-module $L_\alpha$ defined in the following way: the underlying vector space is given by a copy of the field $k$ for each crossing between $T$ and $\alpha$, and the action of an arrow $a \in Q(T)_1$ is the identity morphism if $\alpha$ crosses both $\tau_{s(a)}$ and $\tau_{t(\alpha)}$, zero otherwise.

If $T$ is $\theta$-invariant, then the algebra $A=kQ(T)/I(T)$ is not symmetric.
To obtain a symmetric algebra, the idea is to apply the involution $F_d$ to the polygon:
\begin{definition}[{\cite[Definition 4.16]{ciliberti2024}}]\label{def_f_d}
    $F_d$ is the operation on $\mathbf{P}_{2n+2}$ which consists of the following three steps in order:
    \begin{itemize}
        \item [1)] cut the polygon along $d$;
        \item [2)] reflect the right part with respect to the axis of symmetry of $d$;
        \item [3)] glue  again the right part along $d$.
    \end{itemize}
   
\end{definition}
Let $\rho$ be the reflection of $\mathbf{P}_{2n+2}$ along $d$. Under the bijection $F_d$, $\theta$-orbits correspond (up to isotopy) to $\rho$-orbits. In particular, diameters correspond to $\rho$-invariant diagonals, while pairs of centrally symmetric diagonals correspond to $\rho$-invariant pairs of diagonals which are not orthogonal to $d$ (cf. \cite[Lemma 4.18]{ciliberti2024}). Therefore, the triangulation $T'$ of $\mathbf{P}_{2n+2}$ in the isotopy class of $F_d(T)$ is $\rho$-invariant and contains the diameter $d$, so $Q(T')$ is a cluster-tilted bound symmetric quiver of type $A_{2n-1}$ with a fixed vertex corresponding to $d$ and no fixed arrows, and $A'=kQ(T')/I(T')$ is a symmetric algebra. Furthermore, $A'$ is gentle.

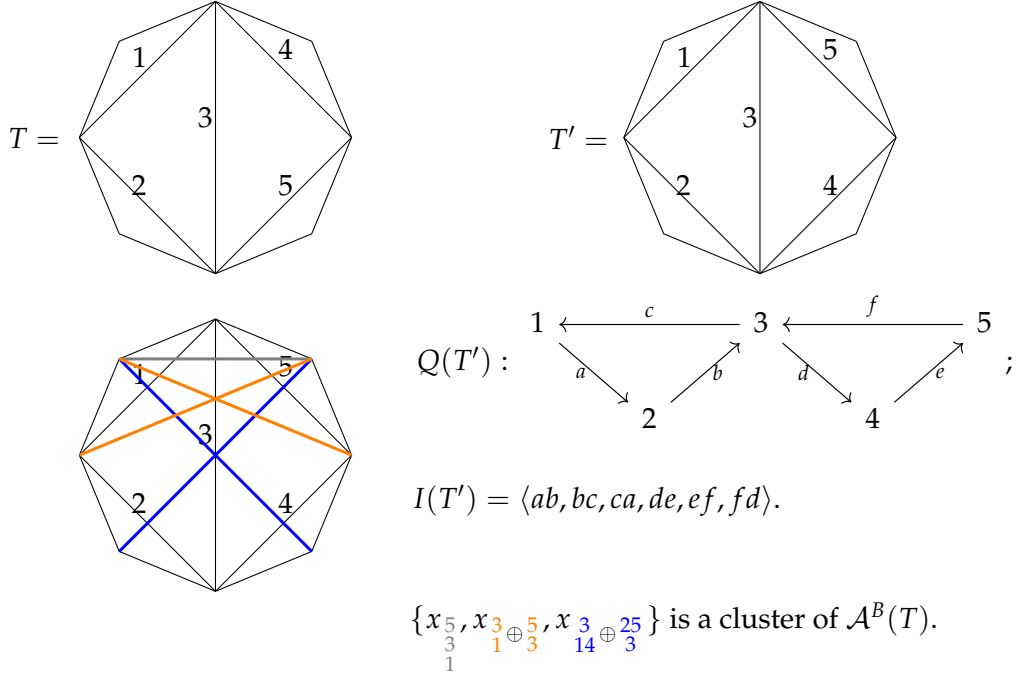
\begin{figure}
        \centering
\begin{tikzpicture}[scale=0.6]
\node at (-4,7) {$T=$};
\begin{scope}[yshift=7cm]
     \draw (90:3cm) -- (135:3cm) -- (180:3cm) -- (225:3cm) -- (270:3cm) -- (315:3cm) -- (360:3cm) -- (45:3cm) -- cycle;
                    \draw (90:3cm) -- node[midway, above, xshift=-1mm, yshift=-1mm] {1} (180:3cm);
                    \draw (270:3cm) -- node[midway, above left,xshift=1.2mm] {2} (180:3cm);
                    \draw (90:3cm) -- node[midway, above left,xshift=1mm] {3} (270:3cm);
                    \draw (90:3cm) -- node[midway, above right, xshift=-2mm] {4} (360:3cm);
                    \draw (270:3cm) -- node[midway, above right, xshift=-2mm] {5} (360:3cm);
                  
\end{scope}
\node at (8,7) {$T'=$};
\begin{scope}[xshift=12cm, yshift=7cm]
      \draw (90:3cm) -- (135:3cm) -- (180:3cm) -- (225:3cm) -- (270:3cm) -- (315:3cm) -- (360:3cm) -- (45:3cm) -- cycle;
                    \draw (90:3cm) -- node[midway, above, xshift=-1mm, yshift=-1mm] {1} (180:3cm);
                    \draw (270:3cm) -- node[midway, above left,xshift=1.2mm] {2} (180:3cm);
                    \draw (90:3cm) -- node[midway, above left,xshift=1mm] {3} (270:3cm);
                    \draw (90:3cm) -- node[midway, above right, xshift=-2mm] {5} (360:3cm);
                    \draw (270:3cm) -- node[midway, above right, xshift=-2mm] {4} (360:3cm);
\end{scope}
    \draw (90:3cm) -- (135:3cm) -- (180:3cm) -- (225:3cm) -- (270:3cm) -- (315:3cm) -- (360:3cm) -- (45:3cm) -- cycle;
                    \draw (90:3cm) -- node[midway, above, xshift=-1mm, yshift=-1mm] {1} (180:3cm);
                    \draw (270:3cm) -- node[midway, above left,xshift=1.2mm] {2} (180:3cm);
                    \draw (90:3cm) -- node[midway, above left,xshift=1mm] {3} (270:3cm);
                    \draw (90:3cm) -- node[midway, above right, xshift=-2mm] {5} (360:3cm);
                    \draw (270:3cm) -- node[midway, above right, xshift=-2mm] {4} (360:3cm);
                    \draw[blue, line width=0.4mm] (135:3cm) -- (315:3cm); 
                    \draw[blue, line width=0.4mm] (225:3cm) -- (45:3cm);
                    \draw[gray, line width=0.4mm] (135:3cm) -- (45:3cm);
                    \draw[orange, line width=0.4mm] (135:3cm) -- (360:3cm);
                    \draw[orange, line width=0.4mm] (45:3cm) -- (180:3cm);
                 
\node at (11,2) {$Q(T'): \begin{tikzcd}
1 \arrow[rd,"a" left] &              & 3 \arrow[ll,"c" above] \arrow[rd,"d" left] &              & 5 \arrow[ll,"f" above] \\
             & 2 \arrow[ru,"b" right] &                         & 4 \arrow[ru,"e" right] &             
\end{tikzcd}$;};
\node at (8.4,-1) {$I(T')=\langle ab,bc, ca,de,ef,fd \rangle$.};
\node at (10.1,-4) {$\{x_{\begin{smallmatrix}
    \textcolor{gray}{5}\\ \textcolor{gray}{3} \\ \textcolor{gray}{1}
\end{smallmatrix}}, x_{\begin{smallmatrix}
    \textcolor{orange}{3}\\ \textcolor{orange}{1} 
\end{smallmatrix}\oplus \begin{smallmatrix}
    \textcolor{orange}{5}\\ \textcolor{orange}{3} 
\end{smallmatrix}},x_{\begin{smallmatrix}
    \textcolor{blue}{3}\\ \textcolor{blue}{14}
\end{smallmatrix}\oplus \begin{smallmatrix}
    \textcolor{blue}{25}\\ \textcolor{blue}{3}
\end{smallmatrix} }\}$ is a cluster of $\mathcal{A}^B (T)$.};

                   
\end{tikzpicture}
        \caption{An example of cluster for a cluster algebra of type $B_3$.}
        \label{fig:ex_typeB}
    \end{figure}

Since $Q(T')$ has no fixed arrows, orthogonal $A'$-modules are either of type I or of type S. Moreover, $\rho$-invariant diagonals $\alpha$ correspond to orthogonal indecomposable $A'$-modules $L_\alpha$ of type I, while $\rho$-invariant pairs $\{\alpha_1,\alpha_2\}$ of diagonals non-orthogonal to $d$ correspond to orthogonal indecomposable $A'$-modules $L_{\alpha_1}\oplus L_{\alpha_2}$ of type S. It follows that non-initial cluster variables of type $B_n$ in $\mathcal{A}_\bullet (T)^B$ correspond to orhogonal indecomposable $A'$-modules. See Figures \ref{fig:ex_typeB} for an example.

Let $\Bar{A}=Q(\Bar{T})/I(\Bar{T})$, where $\Bar{T}=\text{Res}(T')=\text{Res}(T)=\{\tau_1, \dots, \tau_n\}$. The restriction on $\theta$-orbits corresponds to the following operation on orthogonal indecomposable $A'$-modules: 
\begin{definition}[\cite{ciliberti2024}]\label{def_res}
\begin{itemize}
    \item [(i)] Let $N=(V_i,\phi_a,\langle \cdot, \cdot \rangle)$ be an orthogonal indecomposable $A'$-module. Then the \emph{restriction} of $N$ is $\Bar{A}$-module $\text{Res}(N)=(\text{Res}(V),\text{Res}(\phi))$, where $\text{Res}(V)_i=V_i$ if $i \leq n$, $\text{Res}(V)_i=0$ otherwise; and $\text{Res}(\phi)_a=\phi_a$ if $a : i \to j$, with $i,j \leq n$, $\text{Res}(\phi)_a=0$ otherwise. 
    \item [(ii)] Let $v \in \mathbb{Z}_{\geq 0}^{2n-1}$. The \emph{restriction} of $v$, denoted by $\text{Res}(v)$, is the vector of the first $n$ coordinates of $v$.
\end{itemize}
\end{definition}

For an orthogonal indecomposable $A'$-module $N$, $F_N$ and $\bold{g}_N$ denote the $F$-polynomial and the $\bold{g}$-vector of the cluster variable of $\mathcal{A}_\bullet^B(T)$ that corresponds to $N$. On the other hand, $F_{\text{Res}(N)}$ and $\bold{g}_{\text{Res}(N)}$ are the $F$-polynomial and the $\bold{g}$-vector of the $\Bar{A}$-module $\text{Res}(N)$, as in Definitions \ref{defi::F-polynomial} and \ref{def_g_vector}. We can now state and prove the main result of this section:
\begin{theorem}\label{cat_interpr}
Let $N$ be an orthogonal indecomposable $A'$-module. Let $D=\text{diag}(1,\dots,1,2)\in \mathbb{Z}^{n\times n}$.
\begin{itemize}
    \item [(i)] If $\text{Res}(N)=(V_i,\phi_a)$ is indecomposable as $\Bar{A}$-module, then 
    \begin{equation}\label{e_1}
        \text{$F_N=F_{\text{Res}(N)}$,}
    \end{equation}
    and
    \begin{equation}\label{e_2}
 \bold{g}_{N}=\begin{cases}
          \text{$D \bold{g}_{\text{Res}(N)}$ \hspace{1.8cm}if $\textbf{dim} V_n =0$;}\\
          \text{$D \bold{g}_{\text{Res}(N)}+\bold{e}_n$ \hspace{1cm}if $\textbf{dim} V_n \neq 0$.}
      \end{cases}
      \end{equation}
\item [(ii)]      Otherwise, $N=L\oplus \nabla L$ with $\text{dim Ext}^1(\nabla L, L)=1$, and there exists a generating extension
    \begin{equation}
        0 \to L \to G_1 \oplus G_2 \to \nabla L \to 0,
    \end{equation}
    where $G_1$ and $G_2$ are orthogonal indecomposable $A'$-modules of type I. Then
    \begin{equation}\label{e_3}
        F_N= F_{\text{Res}(N)} - \bold{y}^{\text{Res}(\textbf{dim}\underline{\nabla L})}F_{\text{Res}(M)},
    \end{equation}
    and
    \begin{equation}\label{e_4}
        \bold{g}_{N}=D(\bold{g}_{\text{Res}(N)}+\bold{e}_n),
    \end{equation}
where $M$ is the $\leq_{\mathrm{Ext}}$-minimum extension in $A'$ between $\nabla L/\underline{\nabla L}$ and $\overline{L}$. Moreover, $M$ is an orthogonal $A'$-module.
 \end{itemize}  
     
\end{theorem}

\begin{proof}
  Since $\text{Res}(T)=\text{Res}(T')$, (i) is simply a rewriting of Theorem \ref{t1} (i). Assume that $\text{Res}(N)$ is not indecomposable as a $kQ(T')$-module. It follows that $N$ corresponds, on one hand, to a $\rho$-invariant pair $\{(a,\rho(\Bar{b})),(\Bar{b},\rho(a))\}$ of diagonals non-orthogonal to $d$ that cross $d$. On the other hand, via $F_d$, it corresponds to a $\theta$-orbit $[a,b]=\{(a,b),(\Bar{b},\Bar{a})\}$ such that each diagonal of $[a,b]$ crosses $d$ (see Figure \ref{fig_22}).

   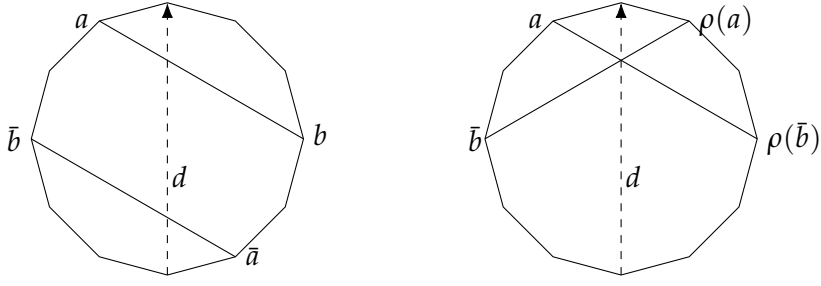
\begin{figure}
        \centering
\begin{tikzpicture}[scale=0.6]
    \draw (90:3cm) -- (120:3cm) -- (150:3cm) -- (180:3cm) -- (210:3cm) -- (240:3cm) -- (270:3cm) -- (300:3cm) -- (330:3cm) -- (360:3cm) -- (30:3cm) -- (60:3cm) --  cycle;
       \draw (120:3cm) -- (360:3cm); 
       \draw (180:3cm) -- (300:3cm); 
         \draw[dashed,-{Latex[length=2mm]}](270:3cm) -- node[midway, left, xshift=0.4cm,yshift=-0.5cm] {$d$} (90:3cm);
         \node at (120:3) [left] {$a$};
         \node at (365:3) [right,yshift=-0.1cm] {$b$};
         \node at (300:3) [right] {$\Bar{a}$};
         \node at (180:3) [left] {$\Bar{b}$};                 
\begin{scope}[xshift=10cm]
   \draw (90:3cm) -- (120:3cm) -- (150:3cm) -- (180:3cm) -- (210:3cm) -- (240:3cm) -- (270:3cm) -- (300:3cm) -- (330:3cm) -- (360:3cm) -- (30:3cm) -- (60:3cm) --  cycle;
     \draw (360:3cm) -- (120:3cm); 
     \draw (60:3cm) -- (180:3cm); 
        \node at (360:3) [right] {$\rho({\Bar{b})}$}; 
       \draw[dashed,-{Latex[length=2mm]}](270:3cm) -- node[midway, left, xshift=0.4cm,yshift=-0.5cm] {$d$} (90:3cm) ;
        \node at (120:3) [left] {$a$};
        \node at (180:3) [left,xshift=0.1cm] {$\Bar{b}$};
         \node at (60:3) [right] {$\rho(a)$};
\end{scope}
                   
\end{tikzpicture}
        \caption{The action of $F_d$ on the $\theta$-orbit $[a,b]$ whose diagonals cross $d$.}
        \label{fig_22}
    \end{figure} 
  
  Let $L=L_{(a,\rho(\Bar{b}))}$, $\nabla L=L_{(\Bar{b},\rho(a))}$, $G_1=L_{(a,\rho(a))}$, $G_2=L_{(\Bar{b},\rho(\Bar{b}))}$. Since $(a,\rho(\Bar{b}))$ and $(\Bar{b},\rho(a))$ cross (in $d$), $[\nabla L, L]^1 = 1$. Moreover, the elementary lamination $L_n$ associated with $\tau_n=d$ crosses both $(a,\rho(a))$ and $(\Bar{b},\rho(\Bar{b}))$, so the middle term of the generating extension between $L$ and $\nabla L$ is $G_1 \oplus G_2$. Therefore, there is a non-split short exact sequence:
\begin{equation}\label{ses}
    0 \to L \to G_1 \oplus G_2 \to \nabla L \to 0,
\end{equation}
where $G_1$ and $G_2$ are $\nabla$-invariant $A'$-modules, and thus orthogonal indecomposable of type I, since $(a,\rho(a))$ and $(\Bar{b},\rho(\Bar{b}))$ are $\rho$-invariant diagonals of $\mathbf{P}_{2n+2}$, and $Q(T')$ has no fixed arrows.

Let $\mathcal{A}_\bullet(T')$ be the cluster algebra of type $A_{2n-1}$ with principal coefficients in $T'$. By Remark \ref{rmk_spec}, we have that 
\begin{align*}
    F_{L\oplus \nabla L}&=F_{G_1\oplus G_2}+\bold{y}^{\textbf{dim}\underline{\nabla L}}F_M,
\end{align*}
where $M$ is the non-trivial extension between $\nabla L/\underline{\nabla L}$ and $\overline{L}$ in $A'$ which is minimal with respect to the $\mathrm{Ext}$-order.
On the other hand, by Proposition \ref{up:skein1},
\begin{align*}
  F_{L\oplus \nabla L}&=F_{G_1\oplus G_2}+\bold{y}^{\bold{d}_{a\rho(a),\Bar{b}\rho(\Bar{b})}}F_{L_{(a,\Bar{b})}\oplus L_{(\rho(a),\rho(\Bar{b}))}}.  
\end{align*}
Thus
\begin{align*}
  \textbf{dim}\underline{\nabla L}&=\bold{d}_{a\rho(a),\Bar{b}\rho(\Bar{b})},   
\end{align*}
and
\begin{align*}
   M&=L_{(a,\Bar{b})}\oplus L_{(\rho(a),\rho(\Bar{b}))}.
\end{align*}
If $\text{Res}([a,b])=\text{Res}(\{(a,\rho(\Bar{b})),(\Bar{b},\rho(a))\})=\{\gamma_1,\gamma_2\}$, then 
\begin{align*}
    &F_{\text{Res}(N)}=F_{L_{\gamma_1}}F_{L_{\gamma_2}}=F_{\gamma_1}F_{\gamma_2},\\
    &\text{Res}(\textbf{dim}\underline{\nabla L})=\text{Res}(\bold{d}_{a\rho(a),\Bar{b}\rho(\Bar{b})})=\bold{d}_{\gamma_1,\gamma_2},\\
    &F_{\text{Res}(M)}=F_{\text{Res}(L_{(a,\Bar{b})}\oplus L_{(\rho(a),\rho(\Bar{b})})}=F_{(a,\Bar{b})}=F_{(a,\theta(b))},\\
    &\bold{g}_{\text{Res}(N)}=\bold{g}_{L_{\gamma_1}}+\bold{g}_{L_{\gamma_2}}=\bold{g}_{\gamma_1}+\bold{g}_{\gamma_2}.
\end{align*}
Therefore, (ii) follows from Theorem \ref{t1} ii). Finally, by Lemma \ref{lemma_M}, $M$ is a $\nabla$-invariant $A'$-module. Since $Q(T')$ has no fixed arrows, $M$ is also orthogonal. 
\end{proof}

\begin{example}
    Let $\mathcal{A}_\bullet^B(T)$ be the cluster algebra of type $B_3$ with principal coefficients in the triangulation $T$ of $\mathbf{P}_8$ in Figure \ref{fig:ex_typeB}. Let $A'=Q(T')/I(T')$ be the corresponding symmetric algebra. We consider the orthogonal indecomposable $A'$-module $N=\begin{smallmatrix}3\\ 14 \end{smallmatrix}\oplus \begin{smallmatrix} 25\\ 3 \end{smallmatrix}$. Let $x_N$ be the cluster variable of $\mathcal{A}_\bullet^B(T)$ that corresponds to $N$, and let $F_N$ and $\bold{g}_N$ denote its $F$-polynomial and its $\bold{g}$-vector respectively. By Theorem \ref{cat_interpr},
    \begin{align*}
        F_N&=F_{\text{Res}(N)}-y_3F_{\text{Res}(\begin{smallmatrix}
            1\\2
        \end{smallmatrix}\oplus \begin{smallmatrix}
            4\\5
        \end{smallmatrix})}=F_{\begin{smallmatrix}
            3\\1
        \end{smallmatrix}\oplus \begin{smallmatrix}
            2\\3
        \end{smallmatrix}}-y_3F_{\begin{smallmatrix}
            1\\2
        \end{smallmatrix}}=1+y_1+2y_1y_3+y_1y_3^2+y_1y_2y_3^2,\\
        \bold{g}_N&= D(\bold{g}_{\text{Res}(N)}+\bold{e}_n)=D(\bold{g}_{\begin{smallmatrix}
            3\\1
        \end{smallmatrix}\oplus \begin{smallmatrix}
            2\\3
        \end{smallmatrix}} + \bold{e}_3)=\begin{pmatrix}
            -1\\0\\0
        \end{pmatrix}.
    \end{align*}
 
\end{example}

\begin{remark}
  The result of Theorem \ref{cat_interpr} is obtained in \cite[Theorem 4.28]{ciliberti2024} for the case where $Q(T')$ does not have oriented cycles, that is, $Q(T')$ is a symmetric orientation of a Dynkin diagram of type $A_{2n-1}$.  
\end{remark}

\section*{Acknowledgments}
The author is grateful to Karin Baur, Giovanni Cerulli Irelli, and Alessandro Contu for valuable discussions and comments.
\addcontentsline{toc}{section}{Acknowledgments}

\printbibliography[heading=bibintoc]

\end{document}